\documentclass[microtype]{gtpart}

\usepackage{pinlabel}

\usepackage[utf8]{inputenc}

\title[Rectification of enriched $\infty$-categories]{Rectification of enriched $\infty$-categories}

\author[R Haugseng]{Rune Haugseng}
\givenname{Rune}
\surname{Haugseng}
\address{Max-Planck-Institut für Mathematik\\
Vivatsgasse 7\\
53111 Bonn\\
Germany}
\email{haugseng@mpim-bonn.mpg.de}
\urladdr{http://people.mpim-bonn.mpg.de/haugseng/}

 	
\arxivreference{1312.3881}
\arxivpassword{4ssjh}

\volumenumber{}
\issuenumber{}
\publicationyear{}
\papernumber{}
\startpage{}
\endpage{}
\doi{}
\MR{}
\Zbl{}
\received{}
\revised{}
\accepted{}
\published{}
\publishedonline{}
\proposed{}
\seconded{}
\corresponding{}
\editor{}
\version{}

\usepackage{eucal}
\usepackage{enumerate}

\newtheorem{thm}{Theorem}[section]
\newtheorem{lemma}[thm]{Lemma}
\newtheorem{propn}[thm]{Proposition}
\newtheorem{cor}[thm]{Corollary}

\theoremstyle{definition}
\newtheorem{defn}[thm]{Definition}
\newtheorem{ex}[thm]{Example}

\newtheorem{remark}[thm]{Remark}

\theoremstyle{remark}

\newcommand{\blank}{\text{--}}

\newcommand{\defterm}[1]{\emph{#1}}
\newcommand{\isoto}{\xrightarrow{\sim}}

\newcommand{\IFF}{if and only if}

\newcommand{\catname}[1]{\ensuremath{\text{\textup{#1}}}}
\newcommand{\txt}[1]{\ensuremath{\text{\textup{#1}}}}
\newcommand{\Set}{\catname{Set}}
\newcommand{\sSet}{\Set_{\Delta}}
\newcommand{\Cat}{\catname{Cat}}
\newcommand{\CatI}{\catname{Cat}_\infty}
\newcommand{\LCatI}{\widehat{\catname{Cat}}_\infty}

\newcommand{\Sp}{\catname{Sp}}

\newcommand{\Fun}{\txt{Fun}}

\newcommand{\Map}{\txt{Map}}
\newcommand{\Hom}{\txt{Hom}}

\newcommand{\op}{\txt{op}}

\newcommand{\icat}{$\infty$-category}
\newcommand{\icats}{$\infty$-categories}
\newcommand{\icatl}{$\infty$-categorical}

\newcommand{\incats}{$(\infty,n)$-categories}

\newcommand{\xto}[1]{\xrightarrow{#1}}

\newcommand{\from}{\leftarrow}

\usepackage{color}

\newcommand{\simp}{\boldsymbol{\Delta}}
\usepackage{tikz}
\usetikzlibrary{matrix,arrows}

\newcommand{\csquare}[8]{ %
\[ %
\begin{tikzpicture} %
\matrix (m) [matrix of math nodes,row sep=3em,column sep=2.5em,text height=1.5ex,text depth=0.25ex] %
{ #1 \pgfmatrixnextcell #2 \\ %
  #3 \pgfmatrixnextcell #4 \\ }; %
\path[->,font=\footnotesize] %
(m-1-1) edge node[auto] {$#5$} (m-1-2)%
(m-1-1) edge node[left] {$#6$} (m-2-1)%
(m-1-2) edge node[auto] {$#7$} (m-2-2)%
(m-2-1) edge node[below] {$#8$} (m-2-2);%
\end{tikzpicture}%
\]%
}

\newcommand{\smallcsquare}[8]{ %
\[ %
\begin{tikzpicture} %
\matrix (m) [matrix of math nodes,row sep=1.5em,column sep=1.25em,text height=1.5ex,text depth=0.25ex] %
{ #1 \pgfmatrixnextcell #2 \\ %
  #3 \pgfmatrixnextcell #4 \\ }; %
\path[->,font=\footnotesize] %
(m-1-1) edge node[auto] {$#5$} (m-1-2)%
(m-1-1) edge node[left] {$#6$} (m-2-1)%
(m-1-2) edge node[auto] {$#7$} (m-2-2)%
(m-2-1) edge node[below] {$#8$} (m-2-2);%
\end{tikzpicture}%
\]%
}

\newcommand{\nolabelcsquare}[4]{\csquare{#1}{#2}{#3}{#4}{}{}{}{}}
\newcommand{\nolabelsmallcsquare}[4]{\smallcsquare{#1}{#2}{#3}{#4}{}{}{}{}}

\newcommand{\liftcsquare}[9]{ %
\[ %
\begin{tikzpicture} %
\matrix (m) [matrix of math nodes,row sep=3em,column sep=2.5em,text height=1.5ex,text depth=0.25ex] %
{ #1 \pgfmatrixnextcell #2 \\ %
  #3 \pgfmatrixnextcell #4 \\ }; %
\path[->,font=\footnotesize] %
(m-1-1) edge node[auto] {$#5$} (m-1-2)%
(m-1-1) edge node[left] {$#6$} (m-2-1)%
(m-1-2) edge node[auto] {$#7$} (m-2-2)%
(m-2-1) edge node[below] {$#8$} (m-2-2);%
\path[->,dashed,font=\footnotesize](m-2-1) edge node[above] {$#9$} (m-1-2);
\end{tikzpicture}%
\]%
}

\newcommand{\ctriangle}[6]{ %
\[ %
\begin{tikzpicture} %
\matrix (m) [matrix of math nodes,row sep=3em,column sep=1.2em,text height=1.5ex,text depth=0.25ex] %
{  \pgfmatrixnextcell #1 \pgfmatrixnextcell \\ %
  #2 \pgfmatrixnextcell \pgfmatrixnextcell #3 \\ }; %
\path[->,font=\footnotesize] %
(m-1-2) edge node[left] {$#4$} (m-2-1)%
(m-1-2) edge node[right] {$#5$} (m-2-3)%
(m-2-1) edge node[below] {$#6$} (m-2-3);%
\end{tikzpicture}%
\]%
}

\newcommand{\opctriangle}[6]{ %
\[ %
\begin{tikzpicture} %
\matrix (m) [matrix of math nodes,row sep=3em,column sep=1.2em,text height=1.5ex,text depth=0.25ex] %
{  #1 \pgfmatrixnextcell \pgfmatrixnextcell #2 \\ %
  \pgfmatrixnextcell #3 \pgfmatrixnextcell \\ %
}; %
\path[->,font=\footnotesize] %
(m-1-1) edge node[above] {$#4$} (m-1-3)%
(m-1-1) edge node[below left] {$#5$} (m-2-2)%
(m-1-3) edge node[below right] {$#6$} (m-2-2);%
\end{tikzpicture}%
\]%
}

\newcommand{\id}{\txt{id}}

\DeclareMathOperator{\colimP}{colim}
\newcommand{\colim}{\mathop{\colimP}}

\DeclareMathOperator{\ob}{ob}

\newcommand{\iopd}{$\infty$-operad}
\newcommand{\iopds}{$\infty$-operads}

\newcommand{\Alg}{\catname{Alg}}

\newcommand{\AlgCat}{\Alg_{\txt{cat}}}

\newcommand{\AlgCatV}{\AlgCat(\mathcal{V})}

\newcommand{\Opd}{\catname{Opd}}
\newcommand{\OpdI}{\Opd_{\infty}}

\newcommand{\OpdIns}{\OpdI^{\txt{ns}}}

\newcommand{\CatIV}{\CatI^{\mathcal{V}}}

\newcommand{\bbGamma}{\boldsymbol{\Gamma}}
\newcommand{\Seg}{\txt{Seg}}

\newcommand{\Mon}{\txt{Mon}}

\newcommand{\RelCat}{\txt{RelCat}}

\begin{document}

\begin{abstract}
  We prove a rectification theorem for enriched \icats{}: If
  $\mathbf{V}$ is a nice monoidal model category, we show that the
  homotopy theory of \icats{} enriched in $\mathbf{V}$ is equivalent
  to the familiar homotopy theory of categories strictly enriched in
  $\mathbf{V}$. It follows, for example, that \icats{} enriched in
  spectra or chain complexes are equivalent to spectral categories and
  dg-categories. A similar method gives a comparison result for
  enriched Segal categories, which implies that the homotopy theories
  of $n$-categories and $(\infty,n)$-categories defined by iterated
  \icatl{} enrichment are equivalent to those of more familiar
  versions of these objects. In the latter case we also include a
  direct comparison with complete $n$-fold Segal spaces. Along the way
  we prove a comparison result for fibrewise simplicial localizations
  potentially of independent use.
\end{abstract}

\begin{asciiabstract}
  We prove a rectification theorem for enriched infinity-categories:
  If V is a nice monoidal model category, we show that the homotopy
  theory of infinify-categories enriched in V is equivalent to the
  familiar homotopy theory of categories strictly enriched in V. It
  follows, for example, that infinity-categories enriched in spectra
  or chain complexes are equivalent to spectral categories and
  dg-categories. A similar method gives a comparison result for
  enriched Segal categories, which implies that the homotopy theories
  of n-categories and (infinity,n)-categories defined by iterated
  infinity-categorical enrichment are equivalent to those of more
  familiar versions of these objects. In the latter case we also
  include a direct comparison with complete n-fold Segal spaces. Along
  the way we prove a comparison result for fibrewise simplicial
  localizations potentially of independent use.
\end{asciiabstract}

\maketitle

\section{Introduction}
In \cite{enr}, David Gepner and I set up a general theory of ``weakly
enriched categories'' --- more precisely, we introduced a notion of
\icats{} \emph{enriched} in a monoidal \icat{}, and constructed an
\icat{} of these objects where the equivalences are the natural
analogue of fully faithful and essentially surjective functors in this
context. In this paper we are interested in the situation where the
monoidal \icat{} we enrich in can be described by a monoidal model
category --- this applies to many, if not most, interesting examples
of monoidal \icats{}. If $\mathbf{V}$ is a model category, then
inverting the weak equivalences $W$ gives an \icat{}
$\mathbf{V}[W^{-1}]$; if $\mathbf{V}$ is a monoidal model category,
then $\mathbf{V}[W^{-1}]$ inherits a monoidal structure, so our theory
produces an \icat{} of $\mathbf{V}[W^{-1}]$-enriched \icats{}. On the
other hand, there is also often a model structure on ordinary
$\mathbf{V}$-enriched categories
(cf. \cite{HTT,BergerMoerdijkEnr,StanculescuEnr,MuroEnr}) where the
weak equivalences are the so-called \emph{DK-equivalences}, namely the
functors that are \emph{weakly fully faithful} (i.e. given by weak
equivalences in $\mathbf{V}$ on morphism objects), and essentially
surjective (up to homotopy). Our main goal in this paper is to prove a
rectification theorem in this setting:
\begin{thm}\label{thm:main} 
  If $\mathbf{V}$ is a nice monoidal model category, then the homotopy
  theory of \icats{} enriched in $\mathbf{V}[W^{-1}]$ is equivalent to
  the homotopy theory of ordinary $\mathbf{V}$-enriched categories
  with respect to the DK-equivalences.
\end{thm}
In particular, $\mathbf{V}[W^{-1}]$-enriched \icats{} can be rectified
to $\mathbf{V}$-categories: every $\mathbf{V}[W^{-1}]$-enriched
\icat{} is equivalent to one coming from a category enriched in $\mathbf{V}$.
We will state and prove a precise version of this result in
\S\ref{sec:enrcatcomp}. The precise meaning of ``nice'' required
applies, for example, to the category of chain complexes over a ring
with the usual projective model structure, and certain model
structures on symmetric spectra. We can therefore conclude that the
\icat{} of spectral categories is equivalent to that of spectral
\icats{}, and the \icat{} of dg-categories to that of \icats{}
enriched in the derived \icat{} of abelian groups.

If $\mathbf{V}$ is a nice Cartesian model category, i.e. a monoidal
model category with respect to the Cartesian product, then the theory
of $\mathbf{V}$-enriched Segal categories, as defined by
Lurie~\cite{LurieGoodwillie} and Simpson~\cite{SimpsonSegCats}, gives
an alternative notion of ``weakly $\mathbf{V}$-enriched
categories''. Using a similar proof strategy we also prove a
comparison result in this setting:
\begin{thm}
  If $\mathbf{V}$ is a nice Cartesian model category, then the
  homotopy theory of \icats{} enriched in $\mathbf{V}[W^{-1}]$ is
  equivalent to the homotopy theory of $\mathbf{V}$-enriched Segal
  categories.
\end{thm}
We will prove a precise version of this theorem in \S\ref{sec:segcatcomp}.
From this we can conclude that the homotopy theories of $n$-categories
and $(\infty,n)$-categories constructed in \cite[\S 6.1]{enr} using
iterated enrichment are equivalent to those constructed as iterated
Segal categories, starting with sets or simplicial sets,
respectively. These are due to Tamsamani and
Pellissier-Hirschowitz-Simpson, and are constructed as model categories
in \cite{SimpsonSegCats}.

Our last main result, which we will prove in \S\ref{sec:compitsegsp},
is a more direct comparison with $(\infty,n)$-categories, generalizing
that between \icats{} enriched in spaces and Segal spaces in \cite[\S
4.4]{enr}:
\begin{thm}
  The homotopy theory of $(\infty,n)$-categories obtained by iterated
  \icatl{} enrichment is equivalent to that of complete $n$-fold Segal
  spaces.
\end{thm}

We now outline the proof of Theorem~\ref{thm:main} and the
organization of the paper. In \cite{enr} we defined enriched \icats{}
in a monoidal \icat{} $\mathcal{V}$ as ``many-object associative
algebras'' in $\mathcal{V}$, or more precisely as algebras for a
``many-object associative operad'' $\simp^{\op}_{X}$, where $X$ is a
space. In \S\ref{sec:review} we briefly review this definition and the
context in which it takes places, namely the theory of non-symmetric \iopds{}.

The first step in the proof of our rectification theorem is to show
that for $X$ a set and $\mathbf{V}$ a nice monoidal model category,
the \icat{} $\Alg_{\simp^{\op}_{X}}(\mathbf{V}[W^{-1}])$ of
$\simp^{\op}_{X}$-algebras in $\mathbf{V}[W^{-1}]$ is equivalent to
the \icat{} obtained by inverting the weakly fully faithful functors
in the category $\Cat_{X}(\mathbf{V})$ of $\mathbf{V}$-categories with
a fixed set of objects $X$. To see this, we first (in
\S\ref{sec:rect}) review Lurie's rectification theorem for associative
algebras (Theorem 4.1.4.4 of \cite{HA}) and observe that it
generalizes to associative algebras in certain non-symmetric monoidal
model categories.

Next, we wish to combine these equivalences to an equivalence of
\icats{} where the sets of objects are allowed to vary. In \cite{enr}
we combined the \icats{} $\Alg_{\simp^{\op}_{X}}(\mathcal{V})$ for all
spaces $X$ to an \icat{} $\AlgCatV$ of \emph{categorical
  algebras}. Here, we consider the \icat{} $\AlgCatV_{\Set}$ of
categorical algebras with \emph{sets} of objects. We will prove that
if $\mathbf{V}$ is a nice monoidal model category, then
$\AlgCat(\mathbf{V}[W^{-1}])_{\Set}$ is equivalent to the \icat{}
obtained from the category $\Cat(\mathbf{V})$ of
$\mathbf{V}$-categories by inverting those morphisms that are weakly
fully faithful and bijective on sets of objects. To see this we need a
technical result about \icatl{} localizations of fibrations of
categories, which we prove in \S\ref{sec:loc}.

The ``correct'' \icat{} of $\mathcal{V}$-\icats{} is not
$\AlgCat(\mathcal{V})$, but rather the \icat{} obtained from this by
inverting the fully faithful and essentially surjective functors. One
of the main results of \cite{enr} was that this is equivalent to the
full subcategory $\CatIV$ of $\AlgCatV$ spanned by those
$\mathcal{V}$-\icats{} that are \emph{complete} in the sense that
their space of objects is equivalent to their classifying space of
equivalences. We also showed, in \cite[Theorem 5.2.17]{enr}, that
inverting the fully faithful and essentially surjective morphisms in
$\AlgCatV$ is equivalent to inverting them in $\AlgCatV_{\Set}$. Since
the DK-equivalences in $\Cat(\mathbf{V})$, if $\mathbf{V}$ is a nice
monoidal model category, correspond to the fully faithful and
essentially surjective functors in
$\AlgCat(\mathbf{V}[W^{-1}])_{\Set}$, we conclude that the \icat{}
obtained from $\Cat(\mathbf{V})$ by inverting the DK-equivalences is
equivalent to $\CatI^{\mathbf{V}[W^{-1}]}$. We will give the details
of the proof we have just sketched in \S\ref{sec:enrcatcomp}, after
the technical preliminaries of \S\ref{sec:rect} and \S\ref{sec:loc}.
We then prove the comparison with Segal categories using a similar
proof in \S\ref{sec:segcatcomp} and the comparison with $n$-fold Segal
spaces in \S\ref{sec:compitsegsp}.

\subsection{Notation}
Much of this paper is based on work of Lurie in \cite{HTT,HA}; we have
generally kept his notation and terminology. In particular, by an
\emph{\icat{}} we mean a quasicategory, i.e. a simplicial set
satisfying certain horn-filling properties. However, in the few cases
where the notation of \cite{enr} differs from that of Lurie we have
kept that of the latter. Here are some hopefully useful reminders:
\begin{itemize}
\item Generic categories are generally denoted by single capital
  bold-face letters ($\mathbf{A},\mathbf{B},\mathbf{C}$) and generic
  \icats{} by single caligraphic letters
  ($\mathcal{A},\mathcal{B},\mathcal{C}$). Specific categories and
  \icats{} both get names in the normal text font: thus the category
  of small $\mathbf{V}$-categories is denoted $\Cat(\mathbf{V})$ and
  the \icat{} of small $\mathcal{V}$-\icats{} is denoted
  $\CatI^{\mathcal{V}}$.
\item $\simp$ is the simplicial indexing category, i.e. the category
  with objects the non-empty ordered sets $[n] = \{0,1,\ldots,n\}$ and
  order-preserving maps as morphisms.
\item A model category is \emph{tractable} if it is combinatorial and
  there exists a set of generating cofibrations that consists of
  morphisms between cofibrant objects.
\item $\sSet$ is the category of simplicial sets, and $\sSet^{+}$ is
  the category of \emph{marked} simplicial sets, i.e. simplicial sets
  equipped with a collection of 1-simplicies including the degenerate
  ones.
\item If $\mathcal{C}$ is an \icat{}, we write $\iota \mathcal{C}$ for
  the \emph{interior} or \emph{underlying space} of $\mathcal{C}$,
  i.e. the largest subspace of $\mathcal{C}$ that is a Kan complex.
\item If $f \colon \mathcal{C} \to \mathcal{D}$ is left adjoint to a
  functor $g \colon \mathcal{D} \to \mathcal{C}$, we will refer to the
  adjunction as $f \dashv g$.
\item $\mathcal{S}$ is the \icat{} of spaces (in the sense of homotopy
  types or $\infty$-groupoids), and $\CatI$ is the \icat{} of
  \icats{}.
\item If $\mathbf{C}$ is a model category, we write
  $\mathbf{C}^{\txt{cof}}$ for the full subcategory of $\mathbf{C}$
  spanned by the cofibrant objects.
\end{itemize}

\subsection{Acknowledgements}
This paper is based on part of my Ph.D. thesis. I thank Haynes Miller,
for being a great Ph.D. advisor in general; Clark Barwick, for several
very helpful conversations, particularly regarding the material in
\S\ref{sec:loc}; Chris Schommer-Pries, for the proof of
Proposition~\ref{propn:SegUnitEq}; and David Gepner, for discussions
of this work and enriched \icats{} in general. I also thank the
American-Scandinavian Foundation and the Norway-America Association
for partially supporting me during the time most of this work was
carried out.

\section{Brief Review of Non-Symmetric $\infty$-Operads and Enriched $\infty$-Categories}\label{sec:review}

To orient the reader, we begin with a brief review of the
non-symmetric version of Lurie's \iopds{} and the definition of
enriched \icats{}. We focus on the essential ideas and do not give
complete technical details of definitions or results; for a more
detailed introduction we refer the reader to \cite[\S 2]{enr}.

The starting point for the theory of non-symmetric \iopds{} is the
\emph{category of operators} of a non-symmetric operad (originally
introduced by May and Thomason for symmetric operads):
\begin{defn}
  Let $\mathbf{O}$ be a coloured non-symmetric operad (or in other words a
  \emph{multicategory}). Its \emph{category of operators}
  $\mathbf{O}^{\otimes}$ has objects (possibly empty) lists
  $(X_{1},\ldots,X_{n})$ of objects of $\mathbf{O}$, and a morphism
  $(X_{1},\ldots,X_{n}) \to (Y_{1},\ldots,Y_{m})$ is given by a morphism
  $\phi \colon [m] \to [n]$ in $\simp$ and for each $i
  = 1,\ldots,m$ a multimorphism in $\mathbf{O}$ from $(X_{\phi(i-1)+1},
  \ldots, X_{\phi(i)})$ to $Y_{i}$.
\end{defn}
There is an obvious
projection $\pi \colon \mathbf{O}^{\otimes} \to \simp^{\op}$, with
the following properties:
\begin{enumerate}[(a)]
\item Recall that a morphism $\phi \colon [n] \to [m]$ in
  $\simp$ is \emph{inert} if it is the inclusion of a subinterval,
  i.e. if $\phi(i) = \phi(0)+i$ for all $i = 0,\ldots,n$. For every
  inert morphism $\phi 
  \colon [n] \to [m]$ and every object $X
  \in \mathbf{O}^{\otimes}$ with $\pi(X) = [m]$, there
  exists a $\pi$-coCartesian morphism $X \to \phi_{!}X$ over $\phi$.
\item Let $\rho_{i} \colon [1] \to [n]$ denote the inert morphism in
  $\simp$ that sends $0$ to $i-1$ and $1$ to $i$. The functors
  $\mathbf{O}^{\otimes}_{[n]} \to \mathbf{O}^{\otimes}_{[1]}$ induced
  by the coCartesian morphisms over $\rho_{i}$ combine to give an
  equivalence of categories
  \[\mathbf{O}^{\otimes}_{[n]} \isoto \prod_{i = 1}^{n} \mathbf{O}^{\otimes}_{[1]}. \]
\item Given objects $X \in \mathbf{O}^{\otimes}_{[n]}$,
  $Y \in \mathbf{O}^{\otimes}_{[m]}$, and a morphism
  $\phi \colon [m] \to [n]$ in $\simp$, the inert maps
  $Y \to \rho_{i,!}Y$ induce an isomorphism
  \[ \Hom_{\mathbf{O}^{\otimes}}^{\phi}(X, Y) \isoto \prod_{i = 1}^{m}
  \Hom_{\mathbf{O}^{\otimes}}^{\rho_{i} \circ \phi}(X, \rho_{i,!}Y),\]
  where $\Hom_{\mathbf{O}^{\otimes}}^{\phi}(X, Y)$ denotes the set of
  morphisms $X \to Y$ in $\mathbf{O}^{\otimes}$ that map to $\phi$ in
  $\simp^{\op}$.
\end{enumerate}
It is not hard to see that these three properties \emph{characterize}
the categories of operators of coloured non-symmetric operads: 
\begin{propn}
  Any functor $\pi \colon \mathbf{C} \to \simp^{\op}$ that satisfies
  (a)--(c) determines a coloured non-symmetric operad that has
  $\mathbf{C}$ as its category of operators. Moreover, under this
  identification morphisms of operads correspond precisely to functors
  over $\simp^{\op}$ that preserve the coCartesian morphisms over the
  inert maps in $\simp^{\op}$.
\end{propn}
Properties (a)--(c) have precise analogues in the theory of \icats{},
and a \emph{non-symmetric \iopd{}} is precisely a functor of
\icats{} $\mathcal{O} \to \simp^{\op}$ with these properties. If
$\mathcal{O}$ and $\mathcal{P}$ are non-symmetric $\infty$-operads in
this sense, it is also easy to define the \icat{}
of \emph{$\mathcal{O}$-algebras} in $\mathcal{P}$:
\begin{defn}
  The \icat{} $\Alg_{\mathcal{O}}(\mathcal{P})$ of
  \emph{$\mathcal{O}$-algebras} in $\mathcal{P}$ is the full
  subcategory of the functor \icat{} $\Fun_{\simp^{\op}}(\mathcal{O},
  \mathcal{P})$ of functors from $\mathcal{O}$ to $\mathcal{P}$ over
  $\simp^{\op}$ spanned by those functors that preserve the
  coCartesian morphisms over inert maps in $\simp^{\op}$.
\end{defn}
The simple definition of the homotopically correct category of
algebras is one of the key advantages of the theory of \iopds{} over
operads enriched in topological spaces or simplicial sets.

An important source of non-symmetric \iopds{} are non-symmetric
operads enriched in simplicial sets or topological spaces: if
$\mathbf{O}$ is a coloured non-symmetric operad enriched in simplicial
sets, all of whose mapping spaces are Kan complexes, then its
simplicial category of operators (defined completely analogously to
the set-based version discussed above) is fibrant, and its coherent
nerve $\mathrm{N}\mathbf{O}^{\otimes} \to \simp^{\op}$ is an \iopd{};
for operads enriched in topological spaces, we simply take the
singular simplicial sets of the mapping spaces first. For example, the
associative operad just gives the identity map $\simp^{\op} \to
\simp^{\op}$, which is easily seen to be equivalent to the \iopd{}
associated to an $A_{\infty}$-operad. This should not be
surprising: in the \icatl{} setting it does not make sense to talk
about ``strict'' associative algebras, the only meaningful notion is
that of an algebra associative up to coherent homotopies, and this
notion is already encoded in algebras for the associative \iopd{}.

We can also recognize monoidal categories from the category of
operators perspective: they are precisely those categories of
operators $\mathbf{C} \to \simp^{\op}$ that are Grothendieck
opfibrations. Analogously we can define a monoidal \icat{} to be a
non-symmetric \iopd{} that is also a coCartesian fibration, but this
can also be reformulated more simply: 
\begin{defn}
  A \emph{monoidal \icat{}} is a coCartesian fibration
  $\mathcal{V}^{\otimes} \to \simp^{\op}$ such that for each $[n] \in
  \simp$ the functor $\mathcal{V}^{\otimes}_{[n]} \to \prod_{i =
    1}^{n} \mathcal{V}^{\otimes}_{[1]}$ induced by the coCartesian
  morphisms over the inert maps $\rho_{i} \colon [1] \to [n]$ is an
  equivalence of \icats{}.
\end{defn}
Using the correspondence between coCartesian fibrations and functors
to the \icat{} $\CatI$ of \icats{}, we get an equivalence between
monoidal \icats{} and \emph{associative monoids} in $\CatI$:
\begin{defn}
  Let $\mathcal{C}$ be an \icat{} with products. An \emph{associative
    monoid} in $\mathcal{C}$ is a functor $\mu \colon \simp^{\op} \to
  \mathcal{C}$ that satisfies the \emph{Segal condition}: for any $[n]
  \in \simp$ the map $\mu([n]) \to \prod_{i = 1}^{n}\mu([1])$ induced
  by the maps $\mu(\rho_{i})$ is an equivalence.
\end{defn}
There is also an equivalence between associative monoids in
$\mathcal{C}$ and algebras for the associative $\infty$-operad in
$\mathcal{C}$ (equipped with the monoidal structure given by the
Cartesian product). In particular, we have:
\begin{propn}
  There are equivalences of \icats{} between associative algebras in
  $\CatI$, associative monoids in $\CatI$, and monoidal \icats{}.
\end{propn}

What we have discussed so far is the non-symmetric variant of
\iopds{}. Lurie's original theory, developed in \cite{HA}, concerns
symmetric \iopds{}. This has a completely analogous motivation, the
only difference is that in the definition of the category of operators
the category $\simp^{\op}$ is replaced by the category $\bbGamma^{\op}$ of pointed
finite sets. In the \icatl{} setting this leads to Lurie's definitions of
symmetric $\infty$-operads and symmetric monoidal \icats{}. As the
non-symmetric theory is the one relevant to the present paper, we
refer the reader to \cite{HA} for more details and do not discuss this
further here.

Instead, we turn to a brief summary of the theory of enriched \icats{}
as introduced in \cite{enr}. Recall that if $\mathbf{V}$ is a monoidal
category, then $\mathbf{V}$-enriched categories with a fixed set $X$
of objects can be regarded as the algebras for a certain non-symmetric
coloured operad $\mathbf{O}_{X}$: 
\begin{defn}
  If $X$ is a set, the multicategory $\mathbf{O}_{X}$ has $X \times X$
  as its set of objects, and the multimorphism sets are defined by \[
  \mathbf{O}_{X}((x_{0}, y_{1}), (x_{1}, y_{2}), \ldots, (x_{n-1},
  y_{n}); (y_{0}, x_{n})) :=
\begin{cases}
  *, & \txt{if $y_{i} = x_{i}$, $i = 0,\ldots, n$,} \\
  \emptyset, & \txt{otherwise.}
\end{cases}
\]
\end{defn}
This suggests that if $\mathcal{V}$ is a monoidal \icat{} then we can
define $\mathcal{V}$-enriched \icats{} with set of objects $X$ to be
algebras in $\mathcal{V}$ for (the non-symmetric \iopd{} associated
to) $\mathbf{O}_{X}$. This is indeed a correct definition, but it
turns out not to be the most convenient to work with --- for instance,
we get a much better-behaved theory of enriched \icats{} if we allow
them to have \emph{spaces} of objects, which is more easily
accomplished with an alternative definition.

We therefore consider \emph{generalized non-symmetric \iopds{}} ---
these are what we obtain by relaxing condition (b) for a category of
operators above to allow $\mathbf{O}^{\otimes}_{[0]}$ to not be just a
point, and instead require $\mathbf{O}^{\otimes}_{[n]}$ to be an
iterated fibre product of $\mathbf{O}^{\otimes}_{[1]}$ over
$\mathbf{O}^{\otimes}_{[0]}$. (The objects that have such categories
of operators in the setting of ordinary categories have been studied
under the names \emph{$\mathbf{fc}$-multicategories} by Leinster and
\emph{virtual double categories} by Cruttwell and Shulman.) For each
set $X$ we can define such a category of operators whose algebras in a
monoidal category (i.e. functors over $\simp^{\op}$ that preserve
coCartesian morphisms over inert maps) are precisely enriched
categories with set of objects $X$:
\begin{defn}
Let $X$ be a set. The category $\simp^{\op}_{X}$ 
has objects lists $(x_{0},\ldots,x_{n})$ of
elements $x_{i} \in X$, and a unique morphism $(x_{0},\ldots,x_{n})
\to (x_{\phi(0)}, \ldots, x_{\phi(m)})$ for each map $\phi \colon [m]
\to [n]$ in $\simp$.
\end{defn}
There is an obvious projection $\simp^{\op}_{X}
\to \simp^{\op}$, and if $\mathbf{V}$ is a monoidal category, then
$\simp^{\op}_{X}$-algebras in the category of operators
$\mathbf{V}^{\otimes}$ are precisely $\mathbf{V}$-enriched categories
with set of objects $X$. This leads to our definition of enriched \icats{}:
\begin{defn}
  If $\mathcal{V}^{\otimes} \to \simp^{\op}$ is a monoidal \icat{},
  then a \emph{$\mathcal{V}$-enriched \icat{}} with set of objects $X$
  is an algebra for the generalized non-symmetric \iopd{}
  $\simp^{\op}_{X}$ in $\mathcal{V}^{\otimes}$.
\end{defn}
The projection $\simp^{\op}_{X} \to \simp^{\op}$ is the Grothendieck
opfibration associated to the functor $\simp^{\op} \to \Set$ that
sends $[n]$ to $X^{\times (n+1)}$ and $\phi \colon [m] \to [n]$ in $\simp$ to
the map $X^{\times(n+1)}\to X^{\times (m+1)}$ that takes $(x_{0},\ldots, x_{n})$ to
$(x_{\phi(0)},\ldots,x_{\phi(m)})$. This has an obvious generalization
where we let $X$ be a space: we simply take the coCartesian fibration
$\simp^{\op}_{X} \to \simp^{\op}$ of the analogous functor
$\simp^{\op} \to \mathcal{X}$ that sends $[n]$ to $X^{\times (n+1)}$.

When $X$ is a set, both $\mathbf{O}_{X}$-algebras and
$\simp^{\op}_{X}$-algebras in a monoidal category $\mathbf{V}$ are
equivalent to $\mathbf{V}$-categories with $X$ as their set of
objects. Similarly, algebras for the non-symmetric \iopd{}
$\mathbf{O}_{X}^{\otimes}$ and the generalized non-symmetric \iopd{}
$\simp^{\op}_{X}$ are equivalent, with the equivalence induced by a
map of generalized \iopds{} (this is a special case of
\cite[Corollary 4.2.8]{enr}):
\begin{propn}\label{propn:NU}
  Suppose $X$ is a set. There is an obvious functor $\nu_{X}$ from
  $\simp^{\op}_{X}$ to $\mathbf{O}_{X}^{\otimes}$ that sends the list
  $(x_{0},\ldots,x_{n})$ to the list $((x_{0},x_{1}),
  (x_{1},x_{2}),\ldots, (x_{n-1}, x_{n}))$. If $\mathcal{V}$ is a
  monoidal \icat{}, then the functor from
  $\Alg_{\mathbf{O}_{X}}(\mathcal{V})$ to
  $\Alg_{\simp^{\op}_{X}}(\mathcal{V})$ given by composition with
  $\nu_{X}$ is an equivalence of \icats{}.
\end{propn}

\section{Rectifying Associative Algebras}\label{sec:rect}
In \cite[\S4.1.4]{HA} Lurie proves a rectification result for
associative algebras: if $\mathbf{V}$ is a nice symmetric monoidal
model category, then the \icat{} of (\icatl{}) associative algebras in
$\mathbf{V}[W^{-1}]$, i.e. the \icat{} of algebras for the
non-symmetric $\infty$-operad $\simp^{\op}$, is equivalent to that
associated to the model category of (strictly) associative algebras in
$\mathbf{V}$, as constructed by Schwede and
Shipley~\cite{SchwedeShipleyAlgMod}. This is proved by showing that
both sides are equivalent to the \icat{} of algebras for the free
associative algebra monad on $\mathbf{V}[W^{-1}]$. In this section we
review this result, and observe that it generalizes slightly to the
setting of non-symmetric monoidal model categories; we will apply this
to enriched categories in \S\ref{sec:enrcatcomp}.


\subsection{Review of Monoidal Model
  Categories}\label{subsec:monmodcat}
In this subsection we briefly review the construction of a monoidal
\icat{} from a monoidal model category; the full details can be found
in \cite[\S 4.1.3]{HA}.

If $\mathbf{V}$ is a simplicial model category, then one way of
constructing an \icat{} from $\mathbf{V}$ is to regard the full
subcategory $\mathbf{V}^{\circ}$ of fibrant-cofibrant objects as a
simplicial category. This is fibrant in the model structure on
simplicial categories, and so its coherent nerve
$\mathrm{N}\mathbf{V}^{\circ}$ is an \icat{}. However, this
construction does not work well with respect to monoidal
structures. We will therefore instead use a more general, but less
explicit, construction, that does not require $\mathbf{V}$ to have a
simplicial enrichment:

\begin{defn}
  Recall that there is a model structure (constructed in
  \cite[\S 3.1.3]{HTT}) on the category $\sSet^{+}$ of marked
  simplicial sets that is Quillen equivalent to the Joyal model
  structure on $\sSet$. In this model category all objects are
  cofibrant and the fibrant objects are precisely those marked
  simplicial sets $(X, S)$ where $X$ is a quasicategory and $S$ is the
  collection of equivalences in $X$. If $\mathcal{C}$ is an \icat{}
  and $W$ is a collection of morphisms in $\mathcal{C}$, then a
  fibrant replacement for the marked simplicial set $(\mathcal{C}, W)$
  in this model structure gives the universal \icat{}
  $\mathcal{C}[W^{-1}]$ obtained from $\mathcal{C}$ by inverting the
  morphisms in $W$.

  If $\mathbf{V}$ is a model category, and $W$ is the class of weak
  equivalences in $\mathbf{V}$, we can therefore define the \icat{}
  $\mathbf{V}[W^{-1}]$ associated to the model category to be a
  fibrant replacement for the marked simplicial set
  $(\mathrm{N}\mathbf{V}, W)$ in this model structure on
  $\sSet^{+}$. Equivalently, we can restrict ourselves to cofibrant,
  fibrant, or fibrant-cofibrant objects and the weak equivalences
  between them. To get monoidal structures on the localization it is
  convenient to consider the cofibrant objects; since this gives an
  \icat{} equivalent to $\mathbf{V}[W^{-1}]$ we will use this notation
  also in this case, despite the slight ambiguity this introduces.
\end{defn}

\begin{defn}
  Let $\mathbf{V}$ be a model category equipped with a biclosed
  monoidal structure. We say that $\mathbf{V}$ is a \emph{monoidal
    model category} if the unit of the monoidal structure is cofibrant
  and the tensor product functor $\otimes \colon \mathbf{V} \times
  \mathbf{V} \to \mathbf{V}$ is a left Quillen bifunctor.
\end{defn}

\begin{remark}\label{rmk:monmodcat}
  Let $\mathbf{V}$ be a model category equipped with a biclosed
  monoidal structure whose unit is cofibrant. If $f
  \colon A \to B$ and $g \colon A' \to B'$ are morphisms in
  $\mathbf{V}$, let $f \square g$ be the induced morphism
  \[ A \otimes B' \amalg_{A \otimes A'} B \otimes A' \to B \otimes
  B';\] this is the \emph{pushout-product} of $f$ and $g$. Then
  $\mathbf{V}$ is a monoidal model category \IFF{} $f \square g$ is a
  cofibration whenever $f$ and $g$ are both cofibrations, and a
  trivial cofibration if either $f$ or $g$ is also a weak equivalence.
\end{remark}

Lurie shows in \cite[Proposition 4.1.3.2]{HA} that the functor that
takes a pair $(\mathcal{C}, W)$ consisting of an \icat{} $\mathcal{C}$
and a collection of morphisms $W$ to the localization
$\mathcal{C}[W^{-1}]$ preserves products. It follows that this functor
preserves $\mathcal{O}$-algebra structures for any $\infty$-operad
$\mathcal{O}$. If $\mathbf{V}$ is a monoidal model category with weak
equivalences $W$, then $(\mathrm{N}\mathbf{V}^{\txt{cof}}, W)$ is an associative
algebra in the \icat{} of such pairs, and so, since a monoidal \icat{}
is the same thing as an algebra for the associative $\infty$-operad in
$\CatI$, we obtain the following key special case of this result:

\begin{propn}[{\cite[Example 4.1.3.6]{HA}}]\label{propn:monmodcatlocmon}
  Let $\mathbf{V}$ be a monoidal model category. Then
  $\mathbf{V}[W^{-1}]$ inherits the structure of a monoidal \icat{}.
\end{propn}

\begin{remark}\label{rmk:unitnotcof}
  The requirement that the unit be cofibrant is often not taken as
  part of the definition of a monoidal model category, as there are
  important examples of model categories with monoidal structures
  where the unit is not cofibrant, but the other requirements for a
  monoidal model category as we have defined it are satisfied. We
  therefore point out that the assumption that $\mathbf{V}$ has a
  cofibrant unit is not essential for
  Proposition~\ref{propn:monmodcatlocmon} to hold. If we drop this
  assumption then $(\mathrm{N}\mathbf{V}, W)$ is still a non-unital
  associative algebra, and so $\mathbf{V}[W^{-1}]$ inherits a
  non-unital monoidal \icat{} structure. It is easy to see that a
  cofibrant replacement for the unit of the monoidal structure in
  $\mathbf{V}$ gives a \emph{quasi-unit} in the sense of
  \cite[Definition 5.4.3.5]{HA} --- roughly speaking, this
  is an object $I$ such that $X \otimes I \simeq X \simeq I \otimes X$
  for every object $X$, but we are not given coherent associativity
  data for combinations of multiple such equivalences. A non-unital
  monoidal \icat{} with a quasi-unit can be extended to a full
  monoidal structure with this as unit in an essentially unique way by
  \cite[Theorem 5.4.3.8]{HA}, and so a monoidal model category without
  a cofibrant unit still induces a monoidal \icat{} structure on its
  associated \icat{}.
\end{remark}

\subsection{Model Categories of Associative Algebras}\label{subsec:rect}
In this subsection we briefly recall the construction of a model
structure on associative algebras, due to Schwede and Shipley, and
observe that it generalizes to non-symmetric monoidal model categories
satisfying an appropriate version of the monoid axiom. First we recall
an observation of Schwede and Shipley on model structures for algebras
over monads:
\begin{defn}
  Let $T$ be a monad on a model category $\mathbf{C}$. We say that $T$
  is an \emph{admissible} monad if there exists a model structure on
  the category $\Alg(T)$ of $T$-algebras where a morphism is a weak
  equivalence or fibration \IFF{} the underlying morphism in
  $\mathbf{C}$ is a weak equivalence or fibration.

  Write $F_{T} \colon \mathbf{C} \rightleftarrows \Alg(T) \colon
  U_{T}$ for the associated adjunction. If $\mathbf{C}$ is a
  combinatorial model category with sets $I$ and $J$ of generating
  cofibrations and trivial cofibrations, we say that $T$ is
  \emph{combinatorially admissible} if it is admissible and the model
  structure on $\Alg(T)$ is combinatorial with $F_{T}(I)$ and
  $F_{T}(J)$ as sets of generating cofibrations and trivial
  cofibrations.
\end{defn}

\begin{remark}
  Given a monad $T$ on $\mathbf{C}$, a model structure on $\Alg(T)$
  where a morphism is a weak equivalence or a fibration \IFF{} its
  underlying morphism in $\mathbf{C}$ is one is unique if it
  exists. Clearly, the existence of such a model structure implies
  certain restrictions on $T$ --- for example, it must preserve
  weak equivalences between cofibrant objects --- but we will not
  attempt to describe these here, as we will only need the following
  admissibility criterion of Schwede and Shipley:  
\end{remark}

\begin{thm}[Schwede-Shipley, {\cite[Lemma 2.3]{SchwedeShipleyAlgMod}}]\label{thm:admcrit}
  Suppose $\mathbf{C}$ is a combinatorial model category and $T$ is a
  filtered-colimit-preserving monad on $\mathbf{C}$, and let $J$ be a
  set of generating trivial cofibrations for $\mathbf{C}$. If the
  underlying morphism in $\mathbf{C}$ of every
  morphism in the weakly saturated class generated by $F_{T}(J)$ in
  $\Alg(T)$ is a weak equivalence, then $T$ is combinatorially
  admissible.
\end{thm}
\begin{remark}\label{rmk:admcrit}
  Since weak equivalences in $\mathbf{C}$ are closed under retracts
  and transfinite composites, the weakly saturated class generated by
  $F_{T}(J)$ will be contained in the weak equivalences provided the
  pushout of any morphism in $F_{T}(J)$ along any morphism in
  $\Alg(T)$ is a weak equivalence.
\end{remark}
In \cite{SchwedeShipleyAlgMod}, Schwede and Shipley analyze such
pushouts in the case of associative algebras. They show that the
pushout is a transfinite composite of pushouts of certain maps, as follows:
\begin{thm}[{Schwede-Shipley~\cite[\S 6]{SchwedeShipleyAlgMod}}]\label{thm:Muro}
  Suppose $\mathbf{C}$ is a combinatorial biclosed monoidal model
  category. Write $\Alg(\mathbf{C})$ for the category of
  associative algebra objects of $\mathbf{C}$ and $F \colon
  \mathbf{C} \rightleftarrows \Alg(\mathbf{C}) : U$ for
  the free algebra functor and forgetful functor. Let $f \colon X \to
  Y$ be a morphism in $\mathbf{C}$ and $g \colon F(X) \to A$ be a
  morphism in $\Alg(\mathbf{C})$. If
  \csquare{F(X)}{F(Y)}{A}{B}{F(f)}{g}{g'}{f'} is a pushout diagram in
  $\Alg(\mathbf{C})$, then there is a sequence of
  morphisms in $\mathbf{C}$
  \[ A = B_{0} \xto{\phi_{1}} B_{1} \xto{\phi_{2}} B_{2 }\cdots \]
  such that $B = \colim_{t} B_{t}$ and $\phi_{t}$ is a pushout of
  $(j \square f)^{\square n} \square j$, where $j$ is the unique
  morphism $\emptyset \to A$.
\end{thm}
Based on this result Schwede and Shipley give a condition --- the
\emph{monoid axiom} --- for the hypothesis of
Theorem~\ref{thm:admcrit} to hold, when the monoidal structure on the
model category $\mathbf{C}$ is \emph{symmetric}, which is true in most
of the interesting examples. However, in the next section we wish to
consider associative algebras in functor cateogories $\Fun(X \times X,
\mathbf{V})$ (where $X$ is a set), equipped with the non-symmetric
``matrix multiplication'' tensor product, for which associative
algebras are precisely $\mathbf{V}$-categories with $X$ as their set
of objects. As noted by Muro~\cite{MuroNSOpds}, the following
non-symmetric version of the monoid axiom applies in this context:
\begin{defn}\label{defn:monoidaxiom}
  Suppose $\mathbf{C}$ is a monoidal model category, and let
  $\mathfrak{U}$ be the set of morphisms in $\mathbf{C}$ of the form
  $f_{1} \square \cdots \square f_{n}$ where each $f_{i}$ is either a
  trivial cofibration or of the form $\emptyset \to X_{i}$ for some
  cofibrant $X_{i}\in \mathbf{C}$, with at least one $f_{i}$ being a trivial
  cofibration. We say that $\mathbf{C}$ satisfies the \emph{monoid
    axiom} if the weakly saturated class $\overline{\mathfrak{U}}$
  generated by $\mathfrak{U}$ is contained in the weak equivalences in
  $\mathbf{C}$.
\end{defn}

\begin{remark}
  Since the pushout-product $(\emptyset \to A) \square f$ is just the
  tensor product $A \otimes f$ for any morphism $f$, the morphisms in
  $\mathfrak{U}$ are all trivial cofibrations in $\mathbf{C}$.
\end{remark}

\begin{remark}
  If $\mathbf{C}$ is \emph{symmetric} monoidal, then we can use the
  symmetry to move all the morphisms of the form $\emptyset \to A$ in
  an element of $\mathfrak{U}$ to one side. Thus, since the pushout
  product of trivial cofibrations in $\mathbf{C}$ is a trivial
  cofibration by Remark~\ref{rmk:monmodcat}, in the symmetric case the
  monoid axiom is equivalent to the corresponding statement where
  $\mathfrak{U}$ consists of morphisms of the form $f \otimes X$ with
  $f$ a trivial cofibration and $X$ a cofibrant object of
  $\mathbf{C}$. This is the original form of the monoid axiom, due to
  Schwede and Shipley.
\end{remark}

\begin{cor}\label{cor:AlgAdm}
  Let $\mathbf{C}$ be a combinatorial biclosed monoidal model category
  that satisfies the monoid axiom. Then the free associative algebra
  monad on $\mathbf{C}$ is combinatorially admissible.
\end{cor}
\begin{proof}
  By Remark~\ref{rmk:admcrit} it suffices to show that if $f \colon X
  \to Y$ is a trivial cofibration in $\mathbf{C}$, $g \colon F(X) \to
  A$ is a morphism in $\Alg(\mathbf{C})$, and
  \csquare{F(X)}{F(Y)}{A}{B}{F(f)}{g}{g'}{f'} is a pushout diagram in
  $\Alg(\mathbf{C})$, then $f'$ is a weak equivalence in
  $\mathbf{C}$. Since $\mathbf{C}$ satisfies the monoid axiom, it
  suffices to show that $f'$ is contained in the weakly saturated
  class $\overline{\mathfrak{U}}$ generated by the class
  $\mathfrak{U}$ from Definition~\ref{defn:monoidaxiom}.

  By Theorem~\ref{thm:Muro}, the morphism $f'$ is a transfinite
  composite of pushouts of morphisms of the form $(j \square
  f)^{\square n} \square j$, where $j$ is the unique morphism
  $\emptyset \to A$, so to show that $f'$ is contained in
  $\overline{\mathfrak{U}}$ it suffices to observe that these
  morphisms are contained in $\mathfrak{U}$ by definition.
\end{proof}

We will also need the following result of Schwede and Shipley:
\begin{cor}\label{cor:AlgCof}
  Let $\mathbf{C}$ be a combinatorial biclosed monoidal model category
  that satisfies the monoid axiom. Then the forgetful functor
  $\Alg(\mathbf{C}) \to \mathbf{C}$ preserves cofibrant objects.
\end{cor}

\subsection{Rectifying Algebras}
We now observe that Lurie's rectification result for associative
algebras also holds for non-symmetric monoidal model categories.
To state the result, we first make the following definition:
\begin{defn}\label{defn:AlgCompFtr}
  Let $\mathbf{C}$ be a left proper tractable biclosed monoidal model
  category that satisfies the monoid axiom. By
  Corollary~\ref{cor:AlgCof}, the forgetful functor from
  associative algebras in $\mathbf{C}$ to $\mathbf{C}$ preserves
  cofibrant objects, so we have a natural functor
  $\Alg(\mathbf{C})^{\txt{cof}} \to
  \Alg(\mathbf{C}^{\txt{cof}}$. It is immediate from the
  construction of the monoidal \icat{} structure on
  $\mathbf{C}[W^{-1}]$ in Proposition~\ref{propn:monmodcatlocmon},
  where $W$ denotes the weak equivalences in $\mathbf{C}$, that there
  is a monoidal functor $\mathbf{C}^{\txt{cof}} \to
  \mathbf{C}[W^{-1}]$, which induces a functor of \icats{}
  $\Alg(\mathbf{C}^{\txt{cof}}) \to
  \Alg_{\simp^{\op}}(\mathbf{C}[W^{-1}])$. The
  composite functor $\Alg(\mathbf{C})^{\txt{cof}} \to
  \Alg_{\simp^{\op}}(\mathbf{C}[W^{-1}])$ clearly
  takes weak equivalences of algebras to equivalences,
  and so induces a functor
  \[\Alg(\mathbf{C})[\hat{W}^{-1}] \to
  \Alg_{\simp^{\op}}(\mathbf{C}[W^{-1}]) \] is an equivalence, where
  $\hat{W}$ denotes the weak equivalences in the model structure on
  $\Alg(\mathbf{C})$.
\end{defn}

\begin{thm}[Lurie]\label{thm:AlgComp}
  Let $\mathbf{C}$ be a left proper tractable biclosed monoidal model
  category that satisfies the monoid axiom. Then the functor of
  \icats{} \[\Alg(\mathbf{C})[\hat{W}^{-1}] \to
  \Alg(\mathbf{C}[W^{-1}]) \] defined
  above is an equivalence.
\end{thm}
The proof is exactly the same as the proof of \cite[Theorem
4.1.4.4]{HA}; in particular, the key technical result \cite[Lemma
4.1.4.13]{HA} generalizes to this context:
\begin{defn}
  Suppose $\mathbf{C}$ is a left proper tractable biclosed monoidal model
  category that satisfies the monoid axiom. Then the forgetful functor $U \colon
  \Alg(\mathbf{C}) \to \mathbf{C}$ takes weak
  equivalences to weak equivalences, by definition of the model
  structure on $\Alg(\mathbf{C})$. The composite functor
  of \icats{} $\Alg(\mathbf{C}) \to \mathbf{C} \to
  \mathbf{C}[W^{-1}]$ thus takes the morphisms in $\hat{W}$ to
  equivalences in $\mathbf{C}[W^{-1}]$ and so factors through a
  unique functor $U_{\infty} \colon
  \Alg(\mathbf{C})[\hat{W}^{-1}] \to \mathbf{C}[W^{-1}]$
  --- this is the functor of \icats{} associated to the right Quillen
  functor $U$.
\end{defn}
\begin{lemma}[Lurie]\label{lem:sifted}
  Suppose $\mathbf{C}$ is a left proper tractable biclosed monoidal model
  category that satisfies the monoid axiom and $\mathbf{I}$ is a small category such
  that $\mathrm{N}\mathbf{I}$ is sifted. Then the forgetful functor
  $U_{\infty} \colon \Alg(\mathbf{C})[\hat{W}^{-1}] \to
  \mathbf{C}[W^{-1}]$ preserves $\mathrm{N}\mathbf{I}$-indexed
  colimits.
\end{lemma}
We omit the proof, as it is exactly the same as that of \cite[Lemma
4.1.4.13]{HA}. We will make use of Lemma~\ref{lem:sifted} in the case
of enriched categories, for which we have the following observation:
\begin{lemma}\label{lem:enrsifted}
  If $\mathbf{V}$ is a left proper tractable biclosed monoidal model
  category satisfying the monoid axiom and $X$ is a set, then there is
  a combinatorial model category structure on the category
  $\Cat_{X}(\mathbf{V})$ such that a morphism is a fibration or weak
  equivalence \IFF{} its image in $\Fun(X \times X, \mathbf{V})$
  is. Moreover, if $\mathbf{I}$ is a small category such that
  $\mathrm{N}\mathbf{I}$ is sifted then the forgetful functor
  \[\Cat_{X}(\mathbf{V})[\txt{FF}_{X}^{-1}] \to \Fun(X \times X,
  \mathbf{V})[W_{X}^{-1}]\] preserves $\mathrm{N}\mathbf{I}$-indexed
  colimits, where $W_{X}$ denotes the class of natural transformations
  that are weak equivalences objectwise.
\end{lemma}
\begin{proof}
  Recall that if $\mathbf{V}$ is a biclosed monoidal category and $X$
  is a set then there is a monoidal structure on $\Fun(X \times X,
  \mathbf{V})$, given by
  \[ (F \otimes G)(x,y) = \coprod_{z \in X} F(x,z) \otimes G(z,y),\]
  such that an associative algebra object in $\Fun(X \times X,
  \mathbf{V})$ is precisely a $\mathbf{V}$-category with objects
  $X$. By \cite[Proposition 10.3]{MuroNSOpds}, if $\mathbf{V}$ is a
  monoidal model category satisfying the monoid axiom, then so is
  $\Fun(X \times X, \mathbf{V})$ equipped with this monoidal
  structure. Applying Corollary~\ref{cor:AlgAdm} and
  Lemma~\ref{lem:sifted} to $\Fun(X \times X, \mathbf{V})$ then
  implies the result.
\end{proof}

\section{Fibrewise Localization}\label{sec:loc}
Suppose we have a functor of ordinary categories $F \colon \mathbf{C}
\to \Cat$ together with a collection $W_{C}$ of weak equivalences in
each category $F(C)$ that is preserved by the functors $F(f)$. Then we
have two ways to construct an \icat{} over $\mathbf{C}$ where these
weak equivalences are inverted: On the one hand we can invert the weak
equivalences in each category $F(C)$ to get a functor $\mathbf{C} \to
\CatI$ that sends $C$ to $F(C)[W_{C}^{-1}]$, which corresponds to a
coCartesian fibration $\mathcal{E} \to \mathbf{C}$. On the other hand,
if $\mathbf{E} \to \mathbf{C}$ is a Grothendieck opfibration
corresponding to $F$ then there is a natural collection $W$ of weak
equivalences in $\mathbf{E}$ induced by those in the fibres, and we
can invert these to get an \icat{} $\mathbf{E}[W^{-1}]$. Our main goal
in this section is to prove that in this situation the natural map
$\mathbf{E}[W^{-1}] \to \mathcal{E}$ is an equivalence of \icats{}.

We will do this in two steps: in \S\ref{subsec:relnerve} we show that
the \icat{} $\mathcal{E}$ here is a fibrant replacement in the
coCartesian model structure on $(\sSet^{+})_{/\mathrm{N}\mathbf{C}}$
for $\mathrm{N}\mathbf{E}$ marked by the edges in $W$, then in \S\ref{subsec:hammock} we use an explicit model
for $\mathbf{E}[W^{-1}]$ to show that this, equipped with a natural
choice of marked edges, is also weakly equivalent to
$(\mathrm{N}\mathbf{E}, W)$.  In addition, we prove in
\S\ref{subsec:topspmodel} that when the weak equivalences in
each category $F(C)$ come from a (combinatorial) model structure, then
there is a (combinatorial) model structure on $\mathbf{E}$ whose weak
equivalences are the morphisms in $W$.

\begin{remark}
  Fibrewise localization has also recently been studied by Hinich in
  \cite{HinichLoc}. His approach is quite different from ours, but
  allows him to prove a comparison analogous to ours also in the more
  general case where the base $\mathbf{C}$ is itself equipped with a
  class of weak equivalences.
\end{remark}

\subsection{The Relative Nerve}\label{subsec:relnerve}
Recall that a \defterm{relative category} is a category $\mathbf{C}$
equipped with a collection of ``weak equivalences'', i.e. a
subcategory $W$ containing all objects and isomorphisms. Write
$\RelCat$ for the obvious category of relative categories;
this has been studied as a model for the theory of
$(\infty,1)$-categories by Barwick and
Kan~\cite{BarwickKanRelCat}. The usual nerve functor from categories
to simplicial sets extends to a functor $L \colon \RelCat \to
\sSet^{+}$ that sends $(\mathbf{C}, W)$ to $(\mathrm{N}\mathbf{C},
\mathrm{N}W_{1})$. In \cite[\S 3.1.3]{HTT} Lurie constructs a model
structure on $\sSet^{+}$ where a fibrant replacement for
$L(\mathbf{C}, W)$ is precisely an \icatl{} localization of $\mathbf{C}$ that
inverts the morphisms in $W$ (marked by the equivalences).

If $\mathbf{C}$ is a category, there is a model structure on
$(\sSet^{+})_{/\mathrm{N}\mathbf{C}}$ where a fibrant object is a
coCartesian fibration marked by its coCartesian morphisms, constructed
in \cite[\S 3.1.3]{HTT}, and in \cite[\S 3.5.2]{HTT} Lurie describes a
right Quillen equivalence $\mathrm{N}^{+}_{\mathbf{C}}$ from the
projective model structure on $\Fun(\mathbf{C}, \sSet^{+})$ to this
model structure on $(\sSet^{+})_{/\mathrm{N}\mathbf{C}}$. Given a
functor $F \colon \mathbf{C} \to \RelCat$ we therefore have two
reasonable ways of constructing a fibrant object of
$(\sSet^{+})_{/\mathrm{N}\mathbf{C}}$:
\begin{enumerate}[(i)]
\item Find a fibrant replacement $\overline{F}$ for the functor $LF
  \colon \mathbf{C} \to \sSet^{+}$, and then form
  $\mathrm{N}^{+}_{\mathbf{C}}\overline{F}$.
\item Construct a Grothendieck opfibration $\mathbf{E} \to \mathbf{C}$
  associated to $F$, regarded as a functor to categories, and write
  $S$ for the collection of 1-simplices in $\mathrm{N}\mathbf{E}$ that
  correspond to composites of (fibrewise) weak equivalences and
  coCartesian morphisms. Then find a fibrant replacement in
  $(\sSet^{+})_{/\mathrm{N}\mathbf{C}}$ for $(\mathrm{N}\mathbf{E}, S)
  \to \mathrm{N}\mathbf{C}$.
\end{enumerate}
Our main goal in this subsection is to prove that these give weakly
equivalent objects. We begin by reviewing the definition of the
functor $\mathrm{N}^{+}_{\mathbf{C}}$:

\begin{defn}
  Let $\mathbf{C}$ be a category. Given a functor $F \colon \mathbf{C}
  \to \sSet$, we define  $\mathrm{N}_{\mathbf{C}}F$ to be the
  simplicial set characterized by the property that a morphism
  $\Delta^{I} \to \mathrm{N}_{\mathbf{C}}F$, where $I$ is a partially
  ordered set, is determined by:
  \begin{enumerate}[(1)]
  \item a functor $\sigma \colon I \to \mathbf{C}$,
  \item for every non-empty subset $J
  \subseteq I$ with maximal element $j$, a map $\tau_{J} \colon
  \Delta^{J} \to F(\sigma(j))$,
  \end{enumerate}
  such that for all subsets $K \subseteq J \subseteq I$ with maximal
  elements $k \in K$
  and $j \in J$, the diagram
  \csquare{\Delta^{K}}{F(\sigma(k))}{\Delta^J}{F(\sigma(j))}{\tau_K}{}{}{\tau_J}
  commutes. This defines a functor $\mathrm{N}_{\mathbf{C}} \colon
  \Fun(\mathbf{C}, \sSet) \to (\sSet)_{/\mathrm{N}\mathbf{C}}$.
\end{defn}

The functor $\mathrm{N}_{\mathbf{C}}$ has a left adjoint, which we
denote \[\mathfrak{F}_{\mathbf{C}} \colon
(\sSet)_{/\mathrm{N}\mathbf{C}} \to \Fun(\mathbf{C}, \sSet).\]

\begin{propn}\label{propn:Fid}
  Let $\pi \colon \mathbf{E} \to \mathbf{C}$ be a functor. Then
  $\mathfrak{F}_{\mathbf{C}}\mathrm{N}\mathbf{E}$ is isomorphic to the
  functor $O_{\pi} \colon \mathbf{C} \to \sSet$ defined by $C \mapsto \mathrm{N}\mathbf{E}_{/C}$.
\end{propn}
\begin{proof}
  We must show that there is a natural isomorphism
  $\Hom(\mathrm{N}\mathbf{E}, \mathrm{N}_{\mathbf{C}}(\blank)) \cong
  \Hom(O_{\pi}, \blank)$; we will do this by defining explicit natural
  transformations \[\phi \colon \Hom(O_{\pi}, \blank) \to
  \Hom(\mathrm{N}\mathbf{E}, \mathrm{N}_{\mathbf{C}}(\blank))\] and
  \[\psi \colon \Hom(\mathrm{N}\mathbf{E},
  \mathrm{N}_{\mathbf{C}}(\blank)) \to \Hom(O_{\pi}, \blank)\] that
  are inverse to each other.

  Given $X \colon \mathbf{C} \to \sSet$ and a natural transformation
  $\eta \colon O_{\pi} \to X$, define $\phi(\eta) \colon
  \mathrm{N}\mathbf{E} \to \mathrm{N}_{\mathbf{C}}X$ to be the
  morphism that sends a simplex $\sigma \colon \Delta^{I} \to
  \mathrm{N}\mathbf{E}$ (which we can identify with a functor $I \to
  \mathbf{E}$) to the simplex of $\mathrm{N}_{\mathbf{C}}X$ determined
  by
  \begin{itemize}
  \item the composite functor $I \to \mathbf{E} \to \mathbf{C}$,
  \item for $J \subseteq I$ with maximal element $j$, the composite
    \[\Delta^{J} \to \mathrm{N}\mathbf{E}_{/\pi(\sigma(j))}
    \xto{\eta_{\pi(\sigma(j))}} X(\pi(\sigma(j))).\]
  \end{itemize}
  
  Conversely, given a map $G \colon \mathrm{N}\mathbf{E} \to
  \mathrm{N}_{\mathbf{C}}X$ of simplicial sets over
  $\mathrm{N}\mathbf{C}$, let $\psi(G)$ be the natural transformation
  $O_{\pi}\to X$ determined as follows: for $C \in \mathbf{C}$, the
  morphism $\psi(G)_{C} \colon \mathrm{N}\mathbf{E}_{/C} \to X(C)$ sends a
  simplex $\sigma \colon \Delta^{I} \to \mathrm{N}\mathbf{E}_{/C}$, where $I$
  has maximal element $i$, to the composite
  \[ \Delta^{I} \xto{\tau} X(\pi \sigma(i)) \xto{X(f)} X(C)\] where
  \begin{itemize}
  \item $\tau$ is the $I$-simplex determined by the image under $G$ of the
    $I$-simplex $\sigma'$ of $\mathrm{N}\mathbf{E}$ underlying $\sigma$,
  \item $f$ is the morphism $\pi(\sigma(i)) \to C$ in
    $\mathbf{C}$ from $\sigma$.
  \end{itemize}
  The remaining data in $G \circ \sigma'$ implies that this defines a
  map of simplicial sets $\mathrm{N}\mathbf{E}_{/C} \to X(C)$, and it
  is also easy to see that $\psi(G)$ is natural in $C$.

  Both $\phi$ and $\psi$ are obviously natural in $X$, and expanding
  out the definitions we see that $\phi\psi = \id$ and $\psi\phi =
  \id$, so we have the required natural isomorphism.
\end{proof}

\begin{defn}
  Let $\mathbf{C}$ be a category. Given a functor $\overline{F} \colon
  \mathbf{C} \to \sSet^{+}$ we define
  $\mathrm{N}_{\mathbf{C}}^{+}\overline{F}$ to be the marked simplicial set
  $(\mathrm{N}_{\mathbf{C}}F, M)$ where $F$ is the underlying functor
  $\mathbf{C} \to \sSet$ of $\overline{F}$, and $M$ is the set of edges
  $\Delta^{1} \to \mathrm{N}_{\mathbf{C}}F$ determined by
  \begin{itemize}
  \item a morphism $f \colon C \to C'$ in $\mathbf{C}$,
  \item a vertex $X \in F(C)$,
  \item a vertex $X' \in F(C')$ and an edge $F(f)(X) \to X'$ that
    is marked in $\overline{F}(C')$.
  \end{itemize}
  This determines a functor $\mathrm{N}_{\mathbf{C}}^{+}\colon
  \Fun(\mathbf{C}, \sSet^{+}) \to
  (\sSet^{+})_{/\mathrm{N}\mathbf{C}}$.
\end{defn}
The functor $\mathrm{N}_{\mathbf{C}}^{+}$ has a left adjoint, which we
denote $\mathfrak{F}_{\mathbf{C}}^{+}$.

\begin{cor}\label{cor:Fplusid}
  Let $\pi \colon \mathbf{E} \to \mathbf{C}$ be a functor, and let $M$
  be a set of edges of $\mathrm{N}\mathbf{E}$ that contains the
  degenerate edges. Then
  $\mathfrak{F}^{+}_{\mathbf{C}}(\mathrm{N}\mathbf{E}, M)$ is
  isomorphic to the functor $\overline{O}_{\pi}$ defined by $C \mapsto
  (\mathrm{N}\mathbf{E}_{/C}, M_{C})$, where $M_{C}$ is the collection
  of edges determined by $E \to E'$ in $\mathbf{E}$ and $\pi(E) \to
  \pi(E') \to C$ in $\mathbf{C}$ such that $\pi(E') \cong C$ and $E \to
  E'$ is in $M$.
\end{cor}
\begin{proof}
  We must show that there is a natural isomorphism
  \[\Hom((\mathrm{N}\mathbf{E}, M), \mathrm{N}^{+}_{\mathbf{C}}(\blank))
  \cong \Hom(\overline{O}_{\pi}, \blank).\] Given $\overline{X} \colon
  \mathbf{C} \to \sSet^{+}$, with underlying functor $X \colon
  \mathbf{C} \to \sSet$, and a morphism $G \colon \mathrm{N}\mathbf{E}
  \to \mathrm{N}_{\mathbf{C}}X$, it is immediate from the definitions
  that $G$ takes an edge $\sigma \colon E \to E'$ of
  $\mathrm{N}\mathbf{E}$ lying over $C \to C'$ in $\mathbf{C}$ to a
  marked edge of $\mathrm{N}_{\mathbf{C}}^{+}\overline{X}$ \IFF{}
  $\phi(G)_{C'}$ takes $\sigma$, regarded as an edge of
  $\mathrm{N}\mathbf{E}_{/C'}$, to a marked edge of
  $\overline{X}(C')$. Thus the natural isomorphism
  $\Hom(\mathrm{N}\mathbf{E}, \mathrm{N}_{\mathbf{C}}X) \cong
  \Hom(O_{\pi}, X)$ of Proposition~\ref{propn:Fid} identifies
  $\Hom((\mathrm{N}\mathbf{E},M),
  \mathrm{N}^{+}_{\mathbf{C}}\overline{X})$, regarded as a subset of
  $\Hom(\mathrm{N}\mathbf{E}, \mathrm{N}_{\mathbf{C}}X)$, with
  $\Hom(\overline{O}_{\pi}, \overline{X})$, regarded as a subset of
  $\Hom(O_{\pi}, X)$.
\end{proof}

\begin{thm}[Lurie, {\cite[Proposition 3.2.5.18]{HTT}}]\ 
  \begin{enumerate}[(i)]
  \item The adjunction $\mathfrak{F}_{\mathbf{C}} \dashv
    \mathrm{N}_{\mathbf{C}}$ is a Quillen equivalence between
    $(\sSet)_{/\mathrm{N}\mathbf{C}}$ equipped with the covariant model
    structure and $\Fun(\mathbf{C}, \sSet)$ equipped with the
    projective model structure.
  \item The adjunction $\mathfrak{F}^{+}_{\mathbf{C}} \dashv
    \mathrm{N}^{+}_{\mathbf{C}}$ is a Quillen equivalence between
    $(\sSet^{+})_{/\mathrm{N}\mathbf{C}}$ equipped with the coCartesian
    model structure and $\Fun(\mathbf{C}, \sSet^{+})$ equipped with
    the projective model structure.
  \end{enumerate}
\end{thm}

\begin{remark}
  By \cite[Lemma 3.2.5.17]{HTT}, the functor
  $\mathfrak{F}^{+}_{\mathbf{C}}$ is naturally weakly equivalent to
  the straightening functor defined in \cite[\S 3.2.1]{HTT}, which
  takes a fibrant functor $\mathbf{C} \to \sSet^{+}$ to the associated
  coCartesian fibration.
\end{remark}

Recall that if $\mathcal{C}$ is an \icat{} we write
$\mathcal{C}^{\natural}$ for the marked simplicial set given by
$\mathcal{C}$ marked by the equivalences, and that if $\mathcal{E} \to
\mathrm{N}\mathbf{C}$ is a coCartesian fibration we write
$\mathcal{E}^{\natural}$ for the object of
$(\sSet^{+})_{/\mathrm{N}\mathbf{C}}$ given by $\mathcal{E}$ marked by
the coCartesian morphisms.
\begin{lemma}\label{lem:NCisopfib}
  Let $F \colon \mathbf{C} \to \Cat$ be a functor. Write $\pi \colon
  \mathbf{E} \to \mathbf{C}$ for the Grothendieck opfibration
  associated to $F$, so that $\mathbf{E}$ has objects pairs $(C \in
  \mathbf{C}, X \in F(C))$ and a morphism $(C,X) \to (D,Y)$ in
  $\mathbf{E}$ is given by a morphism $f \colon C \to D$ in
  $\mathbf{C}$ and a morphism $F(f)(X) \to Y$ in $F(D)$. Then:
  \begin{enumerate}[(i)]
  \item $\mathrm{N}_{\mathbf{C}}(\mathrm{N}F) \to
    \mathrm{N}\mathbf{C}$ is isomorphic to $\mathrm{N}\pi$.
  \item $\mathrm{N}^{+}_{\mathbf{C}}(\mathrm{N}F^{\natural}) \to
    \mathrm{N}\mathbf{C}$ is isomorphic to
    $(\mathrm{N}\mathbf{E})^{\natural} \to \mathrm{N}\mathbf{C}$.
  \end{enumerate}
\end{lemma}
\begin{proof}
  It is clear from the definition of $\mathrm{N}_{\mathbf{C}}$ that
  there is a natural isomorphism between $n$-simplices of
  $\mathrm{N}_{\mathbf{C}}(\mathrm{N}F)$ and $n$-simplices of
  $\mathrm{N}\mathbf{E}$, which proves (i). By definition, an edge of
  $\mathrm{N}^{+}_{\mathbf{C}}(\mathrm{N}F^{\natural})$ is marked if
  it is given by $f \colon C \to C'$ in $\mathbf{C}$, $X \in F(C)$,
  and $F(f)(X) \to X'$ an isomorphism in $F(C')$. Under the
  identification with edges of $\mathrm{N}\mathbf{E}$, such edges
  precisely correspond to the coCartesian edges. This proves (ii).
\end{proof}

\begin{propn}\label{propn:FNLFeq}
  Given $F \colon \mathbf{C} \to \txt{RelCat}$, the counit map
  $\mathfrak{F}_{\mathbf{C}}^{+}\mathrm{N}_{\mathbf{C}}^{+}LF \to LF$
  is a weak equivalence in $\Fun(\mathbf{C}, \sSet^{+})$.
\end{propn}
\begin{proof}
  Since $\Fun(\mathbf{C}, \sSet^{+})$ is equipped with the projective
  model structure, it suffices to show that for all $C \in \mathbf{C}$
  the morphism
  $\mathfrak{F}_{\mathbf{C}}^{+}\mathrm{N}_{\mathbf{C}}^{+}LF(C) \to
  LF(C)$ is a weak equivalence in $\sSet^{+}$. Let $F_{0}$ be the
  underlying functor $\mathbf{C} \to \Cat$, and let $\mathbf{E} \to
  \mathbf{C}$ be the canonical Grothendieck opfibration associated to
  $F_{0}$. Then by Lemma~\ref{lem:NCisopfib} we can identify
  $\mathrm{N}_{\mathbf{C}}^{+}\mathrm{N}F_{0}^{\natural}$ with
  $\mathrm{N}\mathbf{E}^{\natural}$, and so by
  Corollary~\ref{cor:Fplusid} we can identify
  $\mathfrak{F}_{\mathbf{C}}^{+}\mathrm{N}_{\mathbf{C}}^{+}\mathrm{N}F_{0}^{\natural}(C)$
  with $\mathrm{N}\mathbf{E}_{/C}$, marked by the set $M_{C}$ of
  coCartesian morphisms $E \to E'$ such that $\pi(E') = C$.

  The adjunction $\mathfrak{F}^{+}_{\mathbf{C}} \dashv
  \mathrm{N}^{+}_{\mathbf{C}}$ is a Quillen equivalence, so since
  $\mathrm{N}F_{0}^{\natural}$ is fibrant and every object of
  $(\sSet^{+})_{/\mathrm{N}\mathbf{C}}$ is cofibrant, the counit
  $\mathfrak{F}_{\mathbf{C}}^{+}\mathrm{N}_{\mathbf{C}}^{+}\mathrm{N}F_{0}^{\natural}
  \to \mathrm{N}F_{0}^{\natural}$ is a weak equivalence in
  $\Fun(\mathbf{C}, \sSet^{+})$. In particular,
  $(\mathrm{N}\mathbf{E}_{/C}, M_{C}) \to
  \mathrm{N}F_{0}(C)^{\natural}$ is a weak equivalence.

  Let $M'_{C}$ be the set of edges of $\mathrm{N}\mathbf{E}_{/C}$
  corresponding to weak equivalences in $F(C)$. Then we have a pushout
  diagram \nolabelcsquare{(\mathrm{N}\mathbf{E}_{/C},
    M_{C})}{\mathrm{N}F_{0}(C)^{\natural}}{(\mathrm{N}\mathbf{E}_{/C},
    M_{C} \cup M'_{C})}{LF(C),} since both vertical maps are pushouts
  along $\coprod_{f \in M'_{C}} \Delta^{1} \hookrightarrow \coprod_{f
    \in M'_{C}} (\Delta^{1})^{\sharp}$. As the model structure on
  $\sSet^{+}$ is left proper, it follows that
  $(\mathrm{N}\mathbf{E}_{/C}, M_{C} \cup M'_{C}) \to LF(C)$ is a weak
  equivalence.
  
  By Corollary~\ref{cor:Fplusid} we can identify
  $\mathfrak{F}_{\mathbf{C}}^{+}\mathrm{N}_{\mathbf{C}}^{+}LF(C)$
  with the simplicial set $\mathrm{N}\mathbf{E}_{/C}$, marked by the set $M''_{C}$ of
  morphisms $E \to E'$ with $\pi(E') = C$ such that given a
  coCartesian factorization $E \to \overline{E} \to E'$ the morphism
  $\overline{E} \to E'$ is a weak equivalence in $LF(C)$. The obvious map $(\mathrm{N}\mathbf{E}_{/C},
  M_{C} \cup M'_{C}) \to
  \mathfrak{F}_{\mathbf{C}}^{+}\mathrm{N}_{\mathbf{C}}^{+}LF(C)$ is
  therefore marked anodyne, since the edges in $M''_{C}$ are precisely the
  composites of edges in $M_{C}$ and $M'_{C}$. In particular this is
  also a weak equivalence, and so by the 2-out-of-3 property the map
  $\mathfrak{F}_{\mathbf{C}}^{+}\mathrm{N}_{\mathbf{C}}^{+}LF(C) \to
  LF(C)$ is a weak equivalence, as required.  
\end{proof}

\begin{cor}
  Given $F \colon \mathbf{C} \to \txt{RelCat}$, let $LF \to \overline{F}$
  be a fibrant replacement in the projective model structure on
  $\Fun(\mathbf{C}, \sSet^{+})$. Then $\mathrm{N}^{+}_{\mathbf{C}}LF \to
  \mathrm{N}^{+}_{\mathbf{C}}\overline{F}$ is a coCartesian equivalence in
  $(\sSet^{+})_{/\mathrm{N}\mathbf{C}}$.
\end{cor}
\begin{proof}
  The adjunction $\mathfrak{F}_{\mathbf{C}}^{+} \dashv
  \mathrm{N}^{+}_{\mathbf{C}}$ is a Quillen equivalence, so since $\overline{F}$ is
  fibrant and every object of $(\sSet^{+})_{/\mathrm{N}\mathbf{C}}$
  is cofibrant, the morphism $\mathrm{N}^{+}_{\mathbf{C}}LF \to
  \mathrm{N}^{+}_{\mathbf{C}}\overline{F}$ is a weak equivalence \IFF{} the
  adjunct morphism
  $\mathfrak{F}_{\mathbf{C}}^{+}\mathrm{N}^{+}_{\mathbf{C}}LF \to
  \overline{F}$ is a weak equivalence. This follows by the 2-out-of-3
  property, since in the commutative diagram
  \opctriangle{\mathfrak{F}_{\mathbf{C}}^{+}\mathrm{N}^{+}_{\mathbf{C}}LF}{LF}{\overline{F}}{}{}{}
  the morphism $LF \to \overline{F}$ is a weak equivalence by assumption,
  and $\mathfrak{F}_{\mathbf{C}}^{+}\mathrm{N}^{+}_{\mathbf{C}}LF \to
  LF$ is a weak equivalence by Proposition~\ref{propn:FNLFeq}.
\end{proof}

Using Lemma~\ref{lem:NCisopfib} we can equivalently state this as:
\begin{cor}
  Given $F \colon \mathbf{C} \to \txt{RelCat}$, suppose $\pi \colon
  \mathbf{E} \to \mathbf{C}$ is a Grothendieck opfibration corresponding to the
  underlying functor $\mathbf{C} \to \Cat$. Let $M$ be the set of
  morphisms $f \colon E \to E'$ in $\mathbf{E}$ such that given a
  coCartesian factorization $E \to \pi(f)_{!}E \to E'$, the morphism
  $\pi(f)_{!}E \to E'$ is a weak equivalence in $F(\pi(E'))$. Then if
  $LF \to \overline{F}$ is a fibrant replacement in $\Fun(\mathbf{C},
  \sSet^{+})$, there is a coCartesian equivalence
  $(\mathrm{N}\mathbf{E}, M) \to \mathrm{N}^{+}_{\mathbf{C}}\overline{F}$.
\end{cor}

\subsection{The Hammock Localization}\label{subsec:hammock}
Consider a functor $F \colon \mathbf{C} \to \txt{RelCat}$, and let
$\pi \colon \mathbf{E} \to \mathbf{C}$ be an opfibration associated to
the underlying functor $\mathbf{C} \to \Cat$. Our main goal in this
subsection is to prove that inverting the collection $W$ of fibrewise
weak equivalences in $\mathbf{E}$ gives a coCartesian fibration
$\mathbf{E}[W^{-1}] \to \mathbf{C}$. As a corollary, we will also see
that $\mathbf{E}[W^{-1}]$ is the total space of the coCartesian
fibration associated to the functor obtained from $F$ by inverting the
weak equivalences in the relative categories $F(C)$. We will prove
this result by analyzing an explicit model for $\mathbf{E}[W^{-1}]$ as
a simplicial category, namely the \emph{hammock localization}. We
begin by recalling the definition of this, specifically the version
defined in \cite[\S 35]{DwyerHirschhornKanSmith}, and its basic
properties:

\begin{defn}
  A \defterm{zig-zag type} $Z = (Z_{+},Z_{-})$ consists of a
  decomposition $\{1,\ldots,n\} = Z_{+} \amalg Z_{-}$. The
  \defterm{zig-zag category} $\mathbf{ZZ}$ is the category with
  objects zig-zag types and morphisms $Z \to Z'$ given by
  order-preserving morphisms $f \colon \{1,\ldots,n\} \to
  \{1,\ldots,n'\}$ such that $f(Z_{+}) \subseteq Z'_{+}$ and $f(Z_{-})
  \subseteq Z'_{-}$. If $Z$ is a zig-zag type, the associated zig-zag
  category $|Z|$ is the category with objects $0, \ldots, n$ and
  \[ |Z|(i,j) =
  \begin{cases}
    *, & \txt{$i \leq j$, $k \in Z_{+}$ for $k = i+1,\ldots,j$},\\
    *, & \txt{$i \geq j$, $k \in Z_{-}$ for $k = j+1, \ldots,i$},\\
    \emptyset, & \txt{otherwise}.
  \end{cases}
  \]
  This clearly gives a functor $| \blank | \colon
  \mathbf{ZZ} \to \Cat$. If $n$ is an odd integer, we abbreviate
  \[\langle n \rangle := (\{2,4,\ldots,n-1\}, \{1, 3, \ldots, n\})\] and
  if $n$ is an even integer we abbreviate \[\langle n \rangle :=
  (\{1,3,\ldots,n-1\}, \{2,4,\ldots,n\}).\]
\end{defn}

\begin{defn}
  Suppose $(\mathbf{C}, W)$ is a relative category. For $x, y \in
  \mathbf{C}$ and $Z \in \mathbf{ZZ}$ we define
  $\mathbf{L}_{W}\mathbf{C}_{Z}(x,y)$ to be the subcategory of
  $\Fun(|Z|, \mathbf{C})$ whose objects are the functors $F \colon |Z|
  \to \mathbf{C}$ such that $F(0) = x$, $F(n) = y$, and $F(i \to
  (i-1))$ is in $W$ for all $i \in Z_{-}$, and whose morphisms are the
  natural transformations $\eta \colon F \to G$ such that $\eta_{0} =
  \id_{x}$, $\eta_{n} = \id_{y}$, and $\eta_{i}$ is in $W$ for all
  $i$. We write $\mathcal{L}_{W}\mathbf{C}_{Z}(x,y) :=
  \mathrm{N}\mathbf{L}_{W}\mathbf{C}_{Z}(x,y)$.

  This construction gives a functor $\mathbf{ZZ}^{\op} \to \Cat$; we
  let $\mathbf{L}_{W}\mathbf{C}(x,y) \to \mathbf{ZZ}$ be the fibration
  associated to it by the Grothendieck construction. Using
  concatenation of zig-zags we get a strict 2-category
  $\mathbf{L}_{W}\mathbf{C}$ with the same objects as $\mathbf{C}$ and
  with mapping categories $\mathbf{L}_{W}\mathbf{C}(x,y)$; taking
  nerves, this gives a simplicial category $\mathcal{L}_{W}\mathbf{C}$
  whose mapping spaces are $\mathcal{L}_{W}\mathbf{C}(x,y) :=
  \mathrm{N}\mathbf{L}_{W}\mathbf{C}(x,y)$. This simplicial category
  is the \emph{hammock localization} of $(\mathbf{C}, W)$.
\end{defn}

\begin{thm}[Dwyer-Kan]\label{thm:locunivpr}
  Let $(\mathbf{C}, W)$ be a relative category. Then:
  \begin{enumerate}[(i)]
  \item The diagram
    \nolabelcsquare{W}{\mathcal{L}_{W}W}{\mathbf{C}}{\mathcal{L}_{W}\mathbf{C}}
    is a homotopy pushout square in simplicial categories.
  \item If $\mathcal{L}_{W}W \to \overline{\mathcal{L}}_{W}W$ is a fibrant
    replacement in simplicial categories, then
    $\mathrm{N}\overline{\mathcal{L}}_{W}W$ is a Kan complex and
    $\mathrm{N}W \to \mathrm{N}\overline{\mathcal{L}}_{W}W$ is a weak
    equivalence of simplicial sets.
  \end{enumerate}
\end{thm}
\begin{proof}
  (i) follows by combining \cite[Proposition
  35.7]{DwyerHirschhornKanSmith}, \cite[Proposition
  2.2]{DwyerKanCalc}, and \cite[\S 4.5]{DwyerKanLoc} (observe that a
  cofibration in the model structure on simplicial categories with a
  fixed set of objects described in \cite[\S 7]{DwyerKanLoc} is also a
  cofibration in the model structure on simplicial categories).

  To prove (ii), we first observe that it follows from \cite[\S
  9.1]{DwyerKanLoc} that $\mathcal{L}_{W}W$ is a simplicial
  groupoid. If $\mathcal{L}_{W}W \to \overline{\mathcal{L}}_{W}W$ is a
  fibrant replacement in simplicial categories, then
  $\mathrm{N}\overline{\mathcal{L}}_{W}W$ is the nerve of a fibrant
  simplicial groupoid, hence a Kan complex by \cite[Theorem
  3.3]{DwyerKanSimplGpd}. Let $\mathfrak{G}$ denote the left adjoint
  to the nerve of simplicial groupoids, as defined in \cite[\S
  3.1]{DwyerKanSimplGpd}; by \cite[Theorem 3.3]{DwyerKanSimplGpd} the
  morphism $\mathrm{N}W \to \mathrm{N}\overline{\mathcal{L}}_{W}W$ is
  a weak equivalence \IFF{} the adjunct $\mathfrak{G}\mathrm{N}W \to
  \overline{\mathcal{L}}_{W}W$ is a weak equivalence of simplicial
  groupoids. This follows from \cite[\S 5.5]{DwyerKanLoc}, since this
  implies that the mapping spaces in both are the appropriate loop
  spaces of $\mathrm{N}W$.
\end{proof}

\begin{cor}
  Let $(\mathbf{C}, W)$ be a relative category. Suppose
  $\mathcal{L}_{W}\mathbf{C} \to \overline{\mathcal{L}}_{W}\mathbf{C}$ is a
  fibrant replacement in the model category of simplicial
  categories. Then \[L(\mathbf{C},W) \to
  \mathrm{N}\overline{\mathcal{L}}_{W}\mathbf{C}^{\natural}\] is a weak
  equivalence in $\sSet^{+}$.
\end{cor}
\begin{proof}
  We must show that for every \icat{} $\mathcal{D}$, the induced map
  \[ \Map_{\sSet^{+}}(\mathrm{N}\overline{\mathcal{L}}_{W}\mathbf{C}^{\natural}, \mathcal{D}^{\natural})
  \to \Map_{\sSet^{+}}(L(\mathbf{C}, W), \mathcal{D}^{\natural})\] 
  is a weak equivalence of simplicial sets. Observe that
  \[\Map_{\sSet^{+}}(L(\mathbf{C}, W), \mathcal{D}^{\natural}) \simeq
  \Map_{\CatI}(\mathrm{N}\mathbf{C}, \mathcal{D})
  \times_{\Map_{\CatI}(\mathrm{N}W, \mathcal{D})}
  \Map_{\CatI}(\mathrm{N}W, \iota\mathcal{D})\] and
  $\Map_{\CatI}(\mathrm{N}W, \iota\mathcal{D}) \simeq
  \Map_{\mathcal{S}}(\overline{\mathrm{N}W}, \iota\mathcal{D}) \simeq
  \Map_{\CatI}(\overline{\mathrm{N}W}, \mathcal{D})$, where
  $\mathrm{N}W \to \overline{\mathrm{N}W}$ denotes a fibrant
  replacement in the usual model structure on simplicial sets, so this
  is equivalent to requiring
  \nolabelcsquare{\mathrm{N}W}{\overline{\mathrm{N}W}}{\mathrm{N}\mathbf{C}}{\mathrm{N}\overline{\mathcal{L}}_{W}\mathbf{C}}
  to be a homotopy pushout square.  Theorem~\ref{thm:locunivpr}(i)
  implies that
  \nolabelcsquare{\mathrm{N}W}{\mathrm{N}\overline{\mathcal{L}}_{W}W}{\mathrm{N}\mathbf{C}}{\mathrm{N}\overline{\mathcal{L}}_{W}\mathbf{C}}
  is a homotopy pushout square, since $\mathrm{N}$ is a right Quillen
  equivalence and all the objects are fibrant. By
  Theorem~\ref{thm:locunivpr}(ii) we also have that $\mathrm{N}W \to
  \mathrm{N}\overline{\mathcal{L}}_{W}W$ is a fibrant replacement in the usual
  model structure on simplicial sets, so the result follows.
\end{proof}

We now fix a functor $F \colon \mathbf{C} \to \txt{RelCat}$, and let
$\pi \colon \mathbf{E} \to \mathbf{C}$ be a Grothendieck opfibration
associated to the underlying functor $\mathbf{C} \to \Cat$. We say a
morphism $\overline{f} \colon X \to Y$ in $\mathbf{E}$ lying over $f \colon
A \to B$ in $\mathbf{C}$ is a \emph{weak equivalence} if $f$ is an
isomorphism and $f_{!}X \to F$ is a weak equivalence in $F(B)$; write
$W$ for the subcategory of $\mathbf{E}$ whose morphisms are the weak
equivalences. Our goal is to show that the nerve of
$\mathcal{L}_{W}\mathbf{E} \to \mathbf{C}$ is (equivalent to) a
coCartesian fibration. To prove this we need a technical hypothesis on
the relative categories $F(C)$:
\begin{defn}
  A relative category $(\mathbf{C}, W)$ satisfies the
  \emph{two-out-of-three property} if given morphisms $r \colon A \to
  B$ and $s \colon B \to C$ such that two out of $r, s, s\circ r$ are
  in $W$, then so is the third.
\end{defn}
\begin{defn}
 We say that a relative category $\overline{\mathbf{C}} = (\mathbf{C}, W)$ is a
 \defterm{partial model category} if $\overline{\mathbf{C}}$ satisfies
 the two-out-of-three property and $\overline{\mathbf{C}}$ admits a
 \emph{three-arrow calculus}, i.e. there exist subcategories $U, V
 \subseteq W$ such that
 \begin{enumerate}[(i)]
 \item for every zig-zag $A' \xleftarrow{u} A \xto{f} B$ in
   $\mathbf{C}$ with $u \in U$, there exists a functorial zig-zag \[A'
   \xto{f'} B' \xleftarrow{u'} B\] with $u' \in U$ such that
   $u'f = f'u$ and $u'$ is an isomorphism if $u$ is,
 \item for every zig-zag $X \xto{g} Y' \xleftarrow{v} Y$ in
   $\mathbf{C}$ with $v \in V$, there exists a functorial zig-zag \[X
   \xleftarrow{v'} X' \xto{g'} Y\] with $v' \in V$ such that
   $gv' = vg'$ and $v'$ is an isomorphism if $v$ is,
 \item every map $w \in W$ admits a functorial factorization
   $w = vu$ with $u \in U$ and $v \in V$.
 \end{enumerate}
\end{defn}

\begin{remark}
  If $\mathbf{M}$ is a model category (with functorial
  factorizations), then the relative category obtained by equipping
  $\mathbf{M}$ with the weak equivalences in the model structure is a
  partial model category. Similarly, the relative categories obtained
  from the full subcategories $\mathbf{M}^{\txt{cof}}$ of cofibrant
  objects, $\mathbf{M}^{\txt{fib}}$ of fibrant objects, and
  $\mathbf{M}^{\circ}$ of fibrant-cofibrant objects together with the
  weak equivalences between these objects are all partial model
  categories.  The term ``partial model category'' is taken from
  \cite{BarwickKanPartModCat}, but we use the more general definition
  of \cite[36.1]{DwyerHirschhornKanSmith} since the more restrictive
  definition of Barwick and Kan does not include what is for us the
  key example, namely $\mathbf{M}^{\txt{cof}}$ for $\mathbf{M}$ a
  model category.
\end{remark}

\begin{thm}[Dwyer-Kan]\label{thm:partmod}
  Suppose $(\mathbf{C}, W)$ is a partial model category. Then for
  every pair of objects $X,Y \in \mathbf{C}$, the morphism
  $\mathcal{L}_{W}\mathbf{C}_{\langle n \rangle}(X,Y) \to
  \mathcal{L}_{W}\mathbf{C}(X,Y)$ is a weak equivalence of simplicial
  sets for all $n \geq 3$.
\end{thm}
\begin{proof}
  For $n = 3$ this is \cite[Proposition 6.2(i)]{DwyerKanCalc}; the
  general case follows similarly.
\end{proof}

\begin{propn}\label{propn:LHcoCarteq}
  Suppose $F \colon \mathbf{C} \to \txt{RelCat}$ is a functor such
  that $F(C)$ is a partial model category for each $C \in \mathbf{C}$.
  Let $\phi \colon A \to B$ be a morphism in $\mathbf{C}$, and let $X$
  and $Y$ be objects of $\mathbf{E}_{A}$ and $\mathbf{E}_{B}$,
  respectively. Write $\mathcal{L}_{W}\mathbf{E}(X,Y)_{\phi}$ for the
  subspace of $\mathcal{L}_{W}\mathbf{E}(X,Y)$ over $\phi$. The
  morphism
    \[ \overline{\phi}^{*} \colon \mathcal{L}_{W}\mathbf{E}_{B}(\phi_{!}X, Y) \to
  \mathcal{L}_{W}\mathbf{E}(X, Y)_{\phi} \]
  given by composition with a coCartesian morphism $\overline{\phi} \colon
  X \to \phi_{!}X$ is a weak equivalence of simplicial sets.
\end{propn}
\begin{proof}
  It is easy to see that $\mathbf{E}$ is also a partial model
  category. It therefore follows from Theorem~\ref{thm:partmod} that the
  maps $\mathcal{L}_{W}\mathbf{E}_{\langle 4 \rangle}(X,Y)_{\phi} \to
  \mathcal{L}_{W}\mathbf{E}(X,Y)_{\phi}$ and
  $\mathcal{L}_{W}(\mathbf{E}_{B})_{\langle 4 \rangle}(\phi_{!}X, Y)
  \to \mathcal{L}_{W}\mathbf{E}_{B}(\phi_{!}X, Y)$ are weak
  equivalences. Since composition with $\overline{\phi}$ gives a
  functor \[\overline{\phi}^{*}\colon \mathbf{L}_{B} :=
  \mathbf{L}_{W}(\mathbf{E}_{B})_{\langle 4 \rangle}(\phi_{!}X, Y) \to
  \mathbf{L}_{W}\mathbf{E}_{\langle 4 \rangle}(X,Y)_{\phi} =:
  \mathbf{L}\] it therefore suffices to prove that this gives a weak
  equivalence upon taking nerves.

  We will prove this in two steps. Let $\mathbf{L}^{1}$ denote the
  full subcategory of $\mathbf{L}$ spanned by objects
  \[ X = X_{0} \xto{f_{1}} X_{1} \xleftarrow{f_{2}} X_{2} \xto{f_{3}}
  X_{3} \xleftarrow{f_{4}} X_{4} = Y\] such that $X_{i} \in
  \mathbf{E}_{B}$ for $i \geq 1$ and $f_{i}$ lies over $\id_{B}$ in
  $\mathbf{C}$ for $i \geq 2$; then $\overline{\phi}^{*}$ factors as
  \[\mathbf{L}_{B} \xto{f} \mathbf{L}^{1}
  \xto{i} \mathbf{L}.\] We will show that each of these functors induces
  a weak equivalence of nerves.

  First we consider $f \colon \mathbf{L}_{B} \to \mathbf{L}^{1}$,
  given by composition with $\overline{\phi}$. Define $q \colon
  \mathbf{L}^{1} \to \mathbf{L}_{B}$ by sending a zig-zag
  \[ X \xto{g} Z \from Z' \to Y' \from Y\]
  in $\mathbf{L}^{1}$ to
  \[ \phi_{!}X \xto{g'} Z \from Z' \to Y' \from Y\]
  where $X \xto{\overline{\phi}} \phi_{!}X \xto{g'} Z$ is the coCartesian
  factorization of $g$ (which exists since the other maps lie over
  $\id_{B}$). Then it is clear that $qf \simeq \id$ and $fq \simeq \id$,
  so $f$ is an equivalence of categories.

  Next we want to define a functor $p \colon \mathbf{L} \to
  \mathbf{L}^{1}$. Given a zig-zag
  \[ X \xto{g} Z' \from Z \xto{h} Y' \from Y\]
  in $\mathbf{L}$, this lies over
  \[ A \to C' \xleftarrow{\gamma} C \to B' \xleftarrow{\beta} B \]
  where $\gamma$ and $\beta$ are isomorphisms, since weak equivalences
  in $\mathbf{E}$ map to isomorphisms in $\mathbf{C}$. Thus the
  coCartesian maps $Z' \to \gamma^{-1}_{!}Z'$ and $B' \to
  \beta^{-1}_{!}B'$ are isomorphisms, and our
  zig-zag is isomorphic to the zig-zag
  \[ X \to \gamma^{-1}_{!}Z' \from Z \to \beta^{-1}_{!}Y' \from Y.\]
  To define $p$ we may therefore assume that $\beta$ and $\gamma$ are
  identities, in which case $p$ sends
  \[ X \xto{f} Z' \from Z \xto{g} Y' \from Y\]
  lying over
  \[ A \xto{\alpha} C \xleftarrow{\id} C \xto{\psi} B
  \xleftarrow{\id} B \]
  to
  \[ X \to \psi_{!}Z' \from \psi_{!}Z \to Y' \from Y\]
  in $\mathbf{L}^{1}$; this is clearly functorial.

  We wish to prove that $p$ gives an inverse to $i$ after taking
  nerves. It is obvious that $p \circ i \simeq \id$,
  so it suffices to show that $i \circ p$ is homotopic to the identity
  after taking nerves. To see this we consider the natural transformation
  $\eta \colon \mathbf{L} \to \Fun([1],
  \mathbf{L}_{W}\mathbf{E}_{\langle 6 \rangle }(x,y)_{\phi})$ that sends
  our zig-zag to the diagram
  \[
   \begin{tikzpicture} 
 \matrix (m) [matrix of math nodes,row sep=3em,column sep=2.5em,text height=1.5ex,text depth=0.25ex] %
 {
   X & Z' & Z & \psi_{!}Z & \psi_{!}Z & Y' & Y \\
   X & Z'  & Z' & \psi_{!}Z' & \psi_{!}Z & Y' & Y,\\
   };
 \path[->,font=\scriptsize] %
 (m-1-1) edge (m-1-2)
 (m-1-3) edge (m-1-2)
 (m-1-3) edge (m-1-4)
 (m-1-5) edge node[above] {$\id$} (m-1-4)
 (m-1-5) edge (m-1-6)
 (m-1-7) edge (m-1-6)
 (m-2-1) edge (m-2-2)
 (m-2-3) edge node[below] {$\id$} (m-2-2)
 (m-2-3) edge (m-2-4)
 (m-2-5) edge (m-2-4)
 (m-2-5) edge (m-2-6)
 (m-2-7) edge (m-2-6)
 (m-1-1) edge node[auto] {$\id$} (m-2-1)
 (m-1-2) edge node[auto] {$\id$} (m-2-2)
 (m-1-3) edge (m-2-3)
 (m-1-4) edge  (m-2-4)
 (m-1-5) edge node[auto] {$\id$} (m-2-5)
 (m-1-6) edge node[auto] {$\id$} (m-2-6)
 (m-1-7) edge node[auto] {$\id$} (m-2-7)
 ;
 \end{tikzpicture}
  \]
  After composing with the inclusion
  $\mathbf{L}_{W}\mathbf{E}_{\langle 6 \rangle }(x,y)_{\phi} \to
  \mathbf{L}_{W}\mathbf{E}(x,y)_{\phi}$ the functor $\eta_{0}$ is
  clearly linked to the inclusion $\mathbf{L} \to
  \mathbf{L}_{W}\mathbf{E}(x,y)_{\phi}$ by a sequence of natural
  transformations, and similarly $\eta_{1}$ is linked to the composite
  of $i \circ p$ with this inclusion. Since natural transformations
  give homotopies of the induced maps between nerves it follows from
  Theorem~\ref{thm:partmod} that the morphism on nerves induced by $i
  \circ p$ is homotopic to the identity. This completes the proof.
\end{proof}

\begin{thm}\label{thm:LEcoCart}
  Suppose $F \colon \mathbf{C} \to \txt{RelCat}$ is a functor such
  that $F(C)$ is a partial model category for each $C \in
  \mathbf{C}$. There is an \icat{} $\mathbf{E}[W^{-1}]$ such that
  $L(\mathbf{E}, W) \to \mathbf{E}[W^{-1}]^{\natural}$ is a weak
  equivalence in $\sSet^{+}$, and $\mathbf{E}[W^{-1}] \to
  \mathrm{N}\mathbf{C}$ is a coCartesian fibration.
\end{thm}
\begin{proof}
  Let $\mathcal{L}_{W}\mathbf{E} \to \overline{\mathcal{L}}_{W}\mathbf{E}
  \to \mathbf{C}$ denote a factorization of $\mathcal{L}_{W}\mathbf{E}
  \to \mathbf{C}$ as a trivial cofibration followed by a fibration in
  the model category of simplicial categories. Then
  $(\mathrm{N}\overline{\mathcal{L}}_{W}\mathbf{E})^{\natural}$ is a
  fibrant replacement for $L(\mathbf{E}, W)$ in $\sSet^{+}$. By
  \cite[Proposition 2.4.4.3]{HTT} to prove that
  $\mathrm{N}\overline{\mathcal{L}}_{W}\mathbf{E} \to \mathrm{N}\mathbf{C}$
  is equivalent to a coCartesian fibration it suffices to show that
  for each morphism $f \colon C \to D$ in $\mathbf{C}$ and each $X$ in
  $\mathbf{E}_{C}$ we have a homotopy pullback square of simplicial
  sets \csquare{\mathcal{L}_{W}\mathbf{E}(f_{!}X, Y)}{\mathcal{L}_{W}\mathbf{E}(X,
    Y)}{\mathbf{C}(D,B)}{\mathbf{C}(C,B)}{\overline{f}^*}{}{}{f^*} for all $B \in \mathbf{C}$
  and $Y \in \mathbf{E}_{B}$, where $\overline{f} \colon X \to f_{!}X$
  denotes a coCartesian morphism in $\mathbf{E}$ over $f$.

  Since the inclusion of a point in a discrete simplicial set is a Kan
  fibration and the model structure on simplicial sets is right
  proper, given $g \colon D \to B$ the fibres at $\{g\}$ and $\{g
  \circ f\}$ in this diagram are homotopy fibres. To see that the
  diagram is a homotopy pullback square it thus suffices to show that
  composition with $\overline{f}$ induces a weak equivalence
  \[ \mathcal{L}_{W}\mathbf{E}(f_{!}X,Y)_{g} \to
  \mathcal{L}_{W}\mathbf{E}(X,Y)_{g f}\] for all $g \colon D \to
  B$. But by Proposition~\ref{propn:LHcoCarteq}, in the commutative
  diagram
  \ctriangle{\mathcal{L}_{W}\mathbf{E}_{B}((gf)_{!}X,Y)}{\mathcal{L}_{W}\mathbf{E}(f_{!}X,Y)_{g}}{\mathcal{L}_{W}\mathbf{E}(X,Y)_{g
      f}}{}{}{} the diagonal morphisms are both weak equivalences,
  hence by the 2-out-of-3 property so is the horizontal morphism.
\end{proof}

\begin{cor}\label{cor:relcat}
  Suppose $F \colon \mathbf{C} \to \txt{RelCat}$ is a functor such
  that $F(C)$ is a partial model category for each $C \in
  \mathbf{C}$. Let $LF \to \overline{F}$ be a fibrant replacement in
  $\Fun(\mathbf{C}, \sSet^{+})$. Then there is a weak equivalence
  $L(\mathbf{E}, W) \to (\mathrm{N}_{\mathbf{C}}\overline{F})^{\natural}$
  in $\sSet^{+}$.
\end{cor}
\begin{proof}
  By Theorem~\ref{thm:LEcoCart}, there exists a coCartesian
  fibration $\mathbf{E}[W^{-1}] \to \mathrm{N}\mathbf{C}$ with a map
  \[\phi \colon L(\mathbf{E}, W) \to \mathbf{E}[W^{-1}]^{\natural}\]
  that is a weak equivalence in $\sSet^{+}$. The map $\phi$ is also a
  weak equivalence when regarded as a morphism in the over-category
  model structure on
  $(\sSet^{+})_{/\mathrm{N}\mathbf{C}^{\natural}}$. Let
  \[ p_{!} \colon (\sSet^{+})_{/\mathrm{N}\mathbf{C}^{\natural}}
  \rightleftarrows (\sSet^{+})_{/\mathrm{N}\mathbf{C}^{\sharp}} \colon
  p^{*} \] be the adjunction where $p_{!}$ is the identity on the
  underlying marked simplicial sets, and $p^{*}$ forgets the marked
  edges that do not lie over isomorphisms in $\mathbf{C}$. If we equip
  $(\sSet^{+})_{/\mathrm{N}\mathbf{C}^{\natural}}$ with the
  over-category model structure and
  $(\sSet^{+})_{/\mathrm{N}\mathbf{C}^{\sharp}}$ with the coCartesian
  model structure, then this is a Quillen adjunction by
  \cite[Proposition B.2.9]{HTT}, since these functors clearly come
  from a map of categorical patterns.  Since all objects in
  $(\sSet^{+})_{/\mathrm{N}\mathbf{C}^{\natural}}$ are cofibrant, the
  functor $p_{!}$ preserves weak equivalences, and so $\phi$ is also a
  weak equivalence when regarded as a morphism of
  $(\sSet^{+})_{/\mathrm{N}\mathbf{C}^{\sharp}}$.


  Let $M'$ be the set of edges of $\mathrm{N}\mathbf{E}$ corresponding
  to coCartesian morphisms in $\mathbf{E}$, and let
  $\mathbf{E}[W^{-1}]^{+}$ denote the marked simplicial set obtained
  from $\mathbf{E}[W^{-1}]^{\natural}$ by also marking the morphisms
  in the image of $M'$. We have a pushout diagram
  \nolabelcsquare{L(\mathbf{E},
    W)}{\mathbf{E}[W^{-1}]^{\natural}}{(\mathrm{N}\mathbf{E},
    \mathrm{N}W_{1} \cup M')}{\mathbf{E}[W^{-1}]^{+},}
  as both vertical maps are pushouts along $\coprod_{f \in M'}
  \Delta^{1} \hookrightarrow \coprod_{f \in M'}
  (\Delta^{1})^{\sharp}$. Since the model structure on
  $(\sSet^{+})_{/\mathrm{N}\mathbf{C}^{\sharp}}$ is left proper, it follows
   that $(\mathrm{N}\mathbf{E}, \mathrm{N}W_{1} \cup M') \to
   \mathbf{E}[W^{-1}]^{+}$ is a weak equivalence.

   Let $\mathbf{E}[W^{-1}]^{*}$ denote $\mathbf{E}[W^{-1}]$, marked by
   the coCartesian morphisms. These are composites of equivalences and
   morphisms in the image of $M'$, so $\mathbf{E}[W^{-1}]^{+} \to
   \mathbf{E}[W^{-1}]^{*}$ is marked anodyne. Moreover,
   it follows as in the proof of Lemma~\ref{lem:NCisopfib} that
   $\mathrm{N}\mathbf{E}$ marked by the composites of morphisms in
   $\mathrm{N}W_{1}$ and $M'$ is precisely
   $\mathrm{N}^{+}_{\mathbf{C}}LF$, so $(\mathrm{N}\mathbf{E},
   \mathrm{N}W_{1} \cup M') \to \mathrm{N}^{+}_{\mathbf{C}}LF$ is also
   marked anodyne. By the 2-out-of-3 property we therefore have a weak
   equivalence $\mathrm{N}^{+}_{\mathbf{C}}LF \to
   \mathbf{E}[W^{-1}]^{*}$. Thus $\mathbf{E}[W^{-1}]^{*}$ and
   $\mathrm{N}^{+}_{\mathbf{C}}\overline{F}$ are both fibrant replacements
   for $\mathrm{N}^{+}_{\mathbf{C}}LF$, and so are linked by a zig-zag
   of weak equivalences between fibrant objects.

   This implies that the underlying \icats{} $\mathbf{E}[W^{-1}]$ and
   $\mathrm{N}_{\mathbf{C}}\overline{F}$ are equivalent, and so by the
   2-out-of-3 property the map $(\mathrm{N}\mathbf{E}, W) \to
  (\mathrm{N}_{\mathbf{C}}\overline{F})^{\natural}$ is a weak equivalence
  in $\sSet^{+}$, as required.
\end{proof}

\subsection{Total Space Model Structures}\label{subsec:topspmodel}
As before we consider a functor $F \colon \mathbf{C} \to \txt{RelCat}$
and let $\mathbf{E} \to \mathbf{C}$ be an opfibration associated to
$F$. Although not strictly necessary for the applications we are
interested in below, in this subsection we show that if the functor
$F$ is obtained from a suitable functor from $\mathbf{C}$ to the
category of combinatorial model categories, then the relative category
structure on $\mathbf{E}$ considered above also comes from a
combinatorial model category.

\begin{defn}
  Let $\txt{ModCat}^{\mathrm{R}}$ be the category of model categories
  and right Quillen functors. A \defterm{right Quillen presheaf} on a
  category $\mathbf{C}$ is a functor $\mathbf{C}^{\op} \to
  \txt{ModCat}^{\mathrm{R}}$. A right Quillen presheaf is
  \defterm{combinatorial} if it factors through the full subcategory
  of combinatorial model categories.
\end{defn}

\begin{defn}
  Suppose $\mathbf{C}$ is a $\kappa$-accessible category. A right
  Quillen presheaf on $\mathbf{C}$ is \emph{$\kappa$-accessible} if
  for each $\kappa$-filtered diagram $i \colon \mathbf{I} \to
  \mathbf{C}$ with colimit $X$, the category $F(X)$ is the limit of
  the categories $F(i(\alpha))$, and the model structure on $F(X)$ is
  induced by those on $F(i(\alpha))$ in the sense that a map $f \colon
  A \to B$ in $F(X)$ is a (trivial) fibration \IFF{}
  $F(g_{\alpha})(f)$ is a (trivial) fibration in $F(i(\alpha))$ for
  all $\alpha \in \mathbf{I}$, where $g_{\alpha}$ is the canonical
  morphism $i(\alpha) \to X$. We say a right Quillen presheaf $F$ on
  an accessible category $\mathbf{C}$ is \defterm{accessible} if there
  exists a cardinal $\kappa$ such that $\mathbf{C}$ and $F$ are
  $\kappa$-accessible.
\end{defn}

\begin{propn}\label{propn:totalspmodstr}
  Suppose $\mathbf{C}$ is a complete and cocomplete category and $F$ is a
  right Quillen presheaf on $\mathbf{C}$. Let $\pi \colon
  \mathbf{E} \to \mathbf{C}$ be the Grothendieck fibration
  corresponding to $F$. Then there exists a model structure on
  $\mathbf{E}$ such that a morphism $\phi \colon X \to Y$ with image
  $f \colon A \to B$ in $\mathbf{C}$ is
  \begin{itemize}
  \item[(W)] a weak equivalence \IFF{} $f$ is an isomorphism in
    $\mathbf{C}$ and the morphism $f_{!}X \to Y$ is a weak
    equivalence in $F(b)$.
  \item[(F)] a fibration \IFF{} $X \to f^{*}Y$ is a fibration in $F(a)$.
  \item[(C)] a cofibration \IFF{} $f_{!}X \to Y$ is a cofibration in
    $F(b)$.
  \end{itemize}
  Moreover, if $\mathbf{C}$ is a presentable category and $F$ is an
  accessible and combinatorial right Quillen presheaf, then this model
  structure on $\mathbf{E}$ is combinatorial.
\end{propn}

\begin{remark}\label{rmk:weqcond}
  If $f \colon A \to B$ is an isomorphism in $\mathbf{C}$, then $f^{*}
  = F(f)$ is an isomorphism of model categories with inverse
  $f_{!}$. Thus if $\phi \colon X \to Y$ is a morphism in $\mathbf{E}$
  such that $f = \pi(\phi)$ is an isomorphism in $\mathbf{C}$, then
  $f_{!}X \to Y$ is a weak equivalence in $\mathbf{E}_{B}$ \IFF{} $X
  \to f^{*}Y$ is a weak equivalence in $\mathbf{E}_{A}$.
\end{remark}

\begin{remark}
  This model category structure is a particular case of that
  constructed by Roig~\cite{Roig} (and corrected by
  Stanculescu~\cite{StanculescuBifib}), though he does not consider
  the combinatorial case. Roig's construction has also recently been
  significantly generalized by Harpaz and
  Prasma~\cite{HarpazPrasmaGroth}. We include a proof for
  completeness.
\end{remark}

\begin{proof}
  Limits in $\mathbf{E}$ are computed by first taking Cartesian
  pullbacks to the fibre over the limit of the projection of the
  diagram to $\mathbf{C}$, and then taking the limit in that
  fibre. Since all the fibres $\mathbf{E}_{B}$ have limits, it is
  therefore clear that $\mathbf{E}$ has limits. Similarly, since each
  functor $\phi^{*}$ for $\phi$ in $\mathbf{C}$ has a left adjoint,
  and each of the fibres $\mathbf{E}_{B}$ has all colimits, it is clear that
  $\mathbf{E}$ has colimits.

  To show that $\mathbf{E}$ is a model category we must now prove that
  the weak equivalences satisfy the 2-out-of-3 property, and the
  cofibrations and trivial fibrations, as well as the trivial
  cofibrations and fibrations, form weak factorization systems. We
  check the 2-out-of-3 property first: Suppose we have morphisms $\overline{f}
  \colon X \to Y$ and $\overline{g} \colon Y \to Z$ in $\mathbf{E}$ lying
  over $f \colon A \to B$ and $g \colon B \to C$ in $\mathbf{C}$.  If
  two out of the three morphisms $\overline{f}$, $\overline{g}$ and
  $\overline{g}\overline{f}$ are weak equivalences, it is clear that $f$ and $g$
  must be isomorphisms. Thus $g_{!}$ is an isomorphism of model
  categories, and $g_{!}f_{!}X \to g_{!}Y$ is a weak equivalence in
  $\mathbf{E}_{C}$ \IFF{} $f_{!}X \to Y$ is a weak equivalence in
  $\mathbf{E}_{B}$. Combining this with the 2-out-of-3 property for
  weak equivalences in $\mathbf{E}_{C}$ gives the 2-out-of-3 property
  for $\mathbf{E}$.

  We now prove that the cofibrations and trivial fibrations form a
  weak factorization system:
  \begin{enumerate}[(1)]
  \item \emph{Any morphism has a factorization as a cofibration followed by
    a trivial fibration:} Given $\overline{f} \colon X \to Y$ in
    $\mathbf{E}$ lying over $f \colon a \to b$ in $\mathbf{C}$, choose
    a factorization $f_{!}X \to Z \to Y$ of $f_{!}X \to Y$ as a
    cofibration followed by a trivial fibration in
    $\mathbf{E}_{b}$. Then by definition $X \to Z$ is a cofibration
    and $Z \to Y$ is a trivial fibration in $\mathbf{E}$.
  \item \emph{A morphism that has the left lifting property with respect to
    all trivial fibrations is a cofibration:} Suppose $\overline{f} \colon X
    \to Y$, lying over $f \colon A \to B$ in $\mathbf{C}$, has the
    left lifting property with respect to all trivial fibrations. Then
    in particular there exists a lift in all diagrams
    \nolabelcsquare{X}{X'}{Y}{Y'} where $X' \to Y'$ is a trivial
    fibration in $\mathbf{E}_{B}$. By the universal property of
    coCartesian morphisms, this clearly implies that $f_{!}X \to Y$
    has the left lifting property with respect to trivial fibrations
    in $\mathbf{E}_{B}$, and so is a cofibration in
    $\mathbf{E}_{B}$. Thus $\overline{f}$ is a cofibration.
  \item \emph{Cofibrations have the left lifting property with respect to
    trivial fibrations:} Suppose $\overline{f} \colon X \to Y$, lying over
    $f \colon A \to B$ in $\mathbf{C}$, is a cofibration, and $\overline{g}
    \colon X' \to Y'$, lying over $g \colon A' \to B'$, is a trivial
    fibration. Given a commutative diagram
    \csquare{X}{X'}{Y}{Y'}{\overline{\alpha}}{\overline{f}}{\overline{g}}{\overline{\beta}}
    lying over \csquare{A}{A'}{B}{B'}{\alpha}{f}{g}{\beta} we must
    show there exists a lift $Y \to X'$. Since $\overline{g}$ is a trivial
    fibration, $g$ is an isomorphism. Pulling back along $g^{-1}$ and
    pushing forward along $g \alpha = \beta f$ and $\beta$ gives a
    diagram
    \[
    \begin{tikzpicture} 
      \matrix (m) [matrix of math nodes,row sep=3em,column
      sep=2.5em,text height=1.5ex,text depth=0.25ex] 
      {
        X & \beta_{!} f_{!}X & (g^{-1})^{*}X' & X' \\
        Y & \beta_{!}Y &  Y' & Y' \\
      };
      \path[->,font=\scriptsize] %
      (m-1-1) edge (m-1-2)
      (m-1-2) edge (m-1-3)
      (m-1-3) edge (m-1-4)
      (m-2-1) edge (m-2-2)
      (m-2-2) edge (m-2-3)
      (m-2-3) edge (m-2-4)
      (m-1-1) edge (m-2-1)
      (m-1-2) edge (m-2-2)
      (m-1-3) edge (m-2-3)
      (m-1-4) edge (m-2-4)
      ;
    \end{tikzpicture}
    \]
    Here $\beta_{!}f_{!}X \to \beta_{!}Y$ is a cofibration in
    $\mathbf{E}_{B'}$ since $f_{!}X \to Y$ is a cofibration in
    $\mathbf{E}_{B}$ and $\beta_{!}$ is a left Quillen functor, and
    $(g^{-1})^{*}X' \to (g^{-1})^{*}g^{*}Y = Y$ is a trivial fibration
    in $\mathbf{E}_{B'}$ since $X \to g^{*}Y$ is a trivial fibration
    in $\mathbf{E}_{Af'}$ and $(g^{-1})^{*}$ is a right Quillen
    functor. Thus there exists a lift $\beta_{!}Y \to (g^{-1})^{*}X'$
    which gives the desired lift $Y \to X'$.
  \item \emph{A morphism that has the right lifting property with respect to
    all cofibrations is a trivial fibration:} Suppose $\overline{g} \colon
    X' \to Y'$, lying over $g \colon A' \to B'$ in $\mathbf{C}$, has
    the right lifting property with respect to all cofibrations. Then
    in particular there exists a lift in all diagrams
    \nolabelcsquare{X}{X'}{Y}{Y'} where $X \to Y$ is a cofibration in
    $\mathbf{E}_{A'}$. By the universal property of Cartesian
    morphisms, this clearly implies that $X' \to g^{*}Y'$ has the
    right lifting property with respect to cofibrations in
    $\mathbf{E}_{A'}$, and so is a trivial fibration in
    $\mathbf{E}_{A'}$. On the other hand, there exists a lift in the
    diagram \nolabelcsquare{X'}{X'}{g_! X'}{Y',} and projecting this
    down to $\mathbf{C}$ we see that $g$ must be an isomorphism. Thus
    $\overline{g}$ is a trivial fibration in $\mathbf{E}$.
  \end{enumerate}
  The proof that trivial cofibrations and fibrations form a weak
  factorization system is dual to that for cofibrations and trivial
  fibrations, so we omit the details. This completes the proof that
  $\mathbf{E}$ is a model category. 

  Now suppose the right Quillen presheaf $F$ is combinatorial and
  accessible.  It follows from \cite[Theorem 5.3.4]{MakkaiPare} that
  the category $\mathbf{E}$ is accessible, and the functor $\pi$ is
  accessible, thus $\mathbf{E}$ is a presentable category since we
  already proved that it has small colimits.
  
  Let $\kappa$ be a cardinal such that $\mathbf{C}$ is
  $\kappa$-accessible and $\mathbf{E}_{X}$ is $\kappa$-accessible for
  each $\kappa$-compact object $X$ in $\mathbf{C}$. For $X \in
  \mathbf{C}$, let $I_{X}$ and $J_{X}$ be sets of generating
  cofibrations and trivial cofibrations for $\mathbf{E}_{X}$. Let $I$
  and $J$ be the unions of $I_{X}$ and $J_{X}$, respectively, over all
  $\kappa$-compact objects $X \in \mathbf{C}$; then $I$ and $J$ are
  sets.

  Suppose a morphism $\overline{f} \colon X \to Y$, lying over
  $f \colon A \to B$ in $\mathbf{C}$, has the right lifting
  property with respect to the morphisms in $J$; 
  then $X \to f^{*}Y$ is a fibration in $\mathbf{E}_{A}$:
  To see this let $\mathbf{K} \to \mathbf{C}$, $\alpha \mapsto
  A_{\alpha}$, be a $\kappa$-filtered diagram of $\kappa$-compact objects with
  colimit $A$, and let $\gamma_{\alpha} \colon A_{\alpha} \to A$ be
  the canonical morphism. Then $\gamma_{\alpha}^{*}X \to
  \gamma_{\alpha}^{*}f^{*}Y$ has the right lifting property with
  respect to a set of generating trivial cofibrations in $\mathbf{E}_{A_{\alpha}}$,
  and hence this is a fibration in $\mathbf{E}_{A_{\alpha}}$. Since
  the right Quillen presheaf $F$ is $\kappa$-accessible, this implies
  that $X \to f^{*}Y$ is a fibration in
  $\mathbf{E}_{A}$. This means $\overline{f}$ is a fibration in
  $\mathbf{E}$, so $J$ is a set of generating trivial
  cofibrations.

  Similarly, if $\overline{f}$ has the right lifting property with respect
  to the morphisms in $I$, then $X \to f^{*}Y$ is a trivial fibration
  in $\mathbf{E}_{A}$. To find a set of generating cofibrations we
  consider also the set $I'$ of morphisms $\emptyset_{\emptyset} \to
  \emptyset_{C}$ and $\emptyset_{C \amalg C} \to \emptyset_{C}$ where
  $C$ is a $\kappa$-compact object of $\mathbf{C}$ and $\emptyset_{C}$
  denotes the initial object of $\mathbf{E}_{C}$. We claim that if
  $\overline{f} \colon X\to Y$ in $\mathbf{E}$, with image $f \colon A \to
  B$ in $\mathbf{C}$, has the right lifting property with respect to
  the morphisms in $I'$, then $f$ is an isomorphism in
  $\mathbf{C}$. To prove this it suffices to show that for every
  object $C \in \mathbf{C}$ the map $f_{*}\colon \Hom_{\mathbf{C}}(C,
  A') \to \Hom_{\mathbf{C}}(C, B')$ induced by composition with $f$ is
  a bijection; since $\mathbf{C}$ is $\kappa$-presentable it is enough
  to prove this for $C$ a $\kappa$-compact object. Since $\overline{f}$ has the
  right lifting property with respect to $\emptyset_{\emptyset} \to
  \emptyset_{C}$ and every morphism $C \to B$ induces a morphism
  $\emptyset_{C} \to Y$, there exists a lift in the diagram
  \liftcsquare{\emptyset}{A}{C}{B}{}{}{f}{}{} for every map $C\to
  B$; this shows that $f_{*}$ is surjective. Moreover, given two
  morphisms $C \to A$ such that the composites $C \to B$ are equal, we
  get a lift in the diagram \liftcsquare{C \amalg
    C}{A}{C}{B}{}{}{f}{}{} since $\overline{f}$ has the right lifting
  property with respect to $\emptyset_{C\amalg C} \to \emptyset_{C}$;
  thus the two morphisms $C \to A$ must be equal and so $f_{*}$ is
  injective. It follows that if a morphism in $\mathbf{E}$ has the
  right lifting property with respect to the union $I \amalg I'$ then it
  is a trivial fibration, so $I \amalg I'$ is a set of generating
  cofibrations for $\mathbf{E}$. Hence $\mathbf{E}$ is a combinatorial
  model category.
\end{proof}

\begin{remark}
  Let $F$ be a right Quillen presheaf on a category $\mathbf{C}$, and
  let $\mathbf{E} \to \mathbf{C}$ be an opfibration
  associated to the underlying functor to categories. Write $G$ for
  the associated ``left Quillen presheaf'' obtained by passing to left
  adjoints, and let $G^{\txt{cof}} \colon \mathbf{C} \to \txt{RelCat}$ be the
  functor to relative categories obtained by restricting to cofibrant
  objects. Then the full subcategory $\mathbf{E}^{\txt{cof}}$ of
  cofibrant objects in $\mathbf{E}$, with the model structure defined
  above, is the total space of the opfibration associated
  to $G^{\txt{cof}}$, and the weak equivalences in
  $\mathbf{E}^{\txt{cof}}$ are precisely those considered above.
\end{remark}

\section{Rectifying Enriched $\infty$-Categories}\label{sec:enrcatcomp}
Our goal in this section is to prove the main result of this paper:
the homotopy theory of categories enriched in a nice monoidal model
category $\mathbf{V}$ (with respect to the DK-equivalences) is
equivalent to the homotopy theory of \icats{} enriched in the monoidal
\icat{} $\mathbf{V}[W^{-1}]$. We will do this in three steps:
\begin{enumerate}[(1)]
\item We first apply the results of \S\ref{sec:rect} to get an
  equivalence between the \icat{} obtained by inverting the weakly
  fully faithful morphisms in the category $\Cat_{X}(\mathbf{V})$ of
  $\mathbf{V}$-categories with a fixed set of objects $X$ and the
  \icat{} $\Alg_{\simp^{\op}_{X}}(\mathbf{V}[W^{-1}])$ of
  $\simp^{\op}_{X}$-algebras.
\item Next, using the results of \S\ref{sec:loc}, we see that this
  induces an equivalence between the \icat{} obtained by inverting
  those morphisms in the category $\Cat(\mathbf{V})$ of
  small $\mathbf{V}$-categories that are weakly fully faithful and
  bijective on objects and the \icat{}
  $\AlgCat(\mathbf{V}[W^{-1}])_{\Set}$ of categorical algebras in
  $\mathbf{V}[W^{-1}]$ whose spaces of objects are sets.
\item Finally, from this we deduce that the \icat{} obtained by
  inverting the DK-equivalences in $\Cat(\mathbf{V})$ is equivalent to
  the \icat{} $\CatI^{\mathbf{V}[W^{-1}]}$ of
  $\mathbf{V}[W^{-1}]$-\icats{}.
\end{enumerate}
For the first step, the map we wish to prove is an equivalence is
defined as follows:
\begin{defn}
  Suppose $\mathbf{V}$ is a left proper tractable biclosed monoidal
  model category satisfying the monoid axiom, and let $X$ be a
  set. The map of generalized \iopds{} $\nu_{X} \colon \simp^{\op}_{X}
  \to \mathbf{O}_{X}^{\otimes}$ defined in Proposition~\ref{propn:NU}
  gives an equivalence
  \[ \Cat_{X}(\mathbf{V}) \simeq \Alg_{\mathbf{O}_{X}}(\mathbf{V})
  \isoto \Alg_{\simp^{\op}_{X}}(\mathbf{V}).\]
  As in Definition~\ref{defn:AlgCompFtr} the monoidal functor
  $\mathbf{V}^{\txt{cof}} \to \mathbf{V}[W^{-1}]$ induces, since the
  forgetful functor $\Cat_{X}(\mathbf{V}) \to \Fun(X \times X,
  \mathbf{V})$ preserves cofibrant objects by
  Corollary~\ref{cor:AlgCof}, a functor
  \[ \Cat_{X}(\mathbf{V})^{\txt{cof}} \to \Alg_{\simp^{\op}_{X}}(\mathbf{V}[W^{-1}]).\]
  Let $\txt{FF}_{X}$ denote the class of morphisms in
  $\Cat_{X}(\mathbf{V})^{\txt{cof}}$ that are weakly fully faithful, i.e. given by
  weak equivalences on all morphism objects. It is clear that these
  are taken to equivalences in
  $\Alg_{\simp^{\op}_{X}}(\mathbf{V}[W^{-1}])$ by this functor, and so
  there is an induced functor
  \[\eta_{X} \colon \Cat_{X}(\mathbf{V})[\txt{FF}_{X}^{-1}] \to
  \Alg_{\simp^{\op}_{X}}(\mathbf{V}[W^{-1}]).\]
  Moreover, it is clear that this is natural in $X$.
\end{defn}

\begin{propn}\label{propn:CatComp}
  Suppose $\mathbf{V}$ is a left proper tractable biclosed monoidal
  model category satisfying the monoid axiom, and let $X$ be a
  set. The natural
  map \[\eta_{X} \colon \Cat_{X}(\mathbf{V})[\txt{FF}_{X}^{-1}] \to
  \Alg_{\simp^{\op}_{X}}(\mathbf{V}[W^{-1}])\] is an equivalence.
\end{propn}

\begin{proof}[Proof of Proposition~\ref{propn:CatComp}]
  We apply \cite[Corollary 4.7.4.16]{HA} as in the proof of
  \cite[Theorem 4.1.4.4]{HA}: We have a commutative diagram
  \opctriangle{\Cat_X(\mathbf{V})[FF_{X}^{-1}]}{\Alg_{\simp^{\op}_{X}}(\mathbf{V}[W^{-1}])}{\Fun(X\times
    X, \mathbf{V}[W^{-1}]),}{\eta_{X}}{U}{V} where $U_{\infty}$ is the
  functor of \icats{} associated to the forgetful functor
  \[U \colon \Cat_{X}(\mathbf{V}) \to \Fun(X \times X, \mathbf{V}),\] which is a
  right Quillen functor, and $V$ is given by
  restricting $\simp^{\op}_{X}$-algebras to the fibre
  $(\simp^{\op}_{X})_{[1]} \simeq X \times X$. Then we observe:
  \begin{enumerate}[(a)]
  \item The \icat{} $\Cat_{X}(\mathbf{V})[\txt{FF}_{X}^{-1}]$ is
    presentable by \cite[Proposition 1.3.4.22]{HA}, and the \icat{}
    $\Alg_{\simp^{\op}_{X}}(\mathbf{V}[W^{-1}])$ is presentable by
    \cite[Corollary B.5.7]{enr} since $\mathbf{V}[W^{-1}]$ is
    presentable by \cite[Proposition 1.3.4.22]{HA} and the induced
    tensor product on $\mathbf{V}[W^{-1}]$ preserves colimits in each
    variable by \cite[Lemma 4.1.4.8]{HA}.
  \item The functor $V$ admits a left adjoint $G$ by
    \cite[Theorem B.4.6]{enr}.
  \item The functor $U_{\infty}$ also admits a left adjoint
    $F_{\infty}$ since it arises from a right Quillen functor.
  \item The functor $V$ is conservative by \cite[Lemma B.5.5]{enr} and
    preserves sifted colimits by \cite[Corollary
    B.5.4]{enr}.
  \item The functor $U_{\infty}$ is conservative by the definition of
    the weak equivalences in $\Alg(\mathbf{V})$, and preserves sifted
    colimits by Lemma~\ref{lem:enrsifted}.
  \item The canonical map $V \circ G \to U_{\infty} \circ
    F_{\infty}$ is an equivalence since both induce, on the level of
    homotopy categories, the free $\mathbf{V}$-category monad \[\Phi
    \mapsto \coprod_{n \geq 0}\coprod_{x_{0},\ldots,x_{n} \in X}
    \Phi(x_{0},x_{1})\otimes \cdots \otimes \Phi(x_{n-1},x_{n}).\]
    This is obvious for $U_{\infty} \circ F_{\infty}$ and for
    $V \circ G$ it follows by
    \cite[Proposition B.4.9]{enr}.
  \end{enumerate}
  The hypotheses of \cite[Corollary 4.7.4.16]{HA} thus hold, which
  implies that the morphism in question is an equivalence.
\end{proof}

For the second step, let us first define the class of maps in
$\Cat(\mathbf{V})$ that we will invert:
\begin{defn}
We say that a functor $F \colon \mathbf{C} \to \mathbf{D}$ of
$\mathbf{V}$-categories is \emph{weakly fully faithful} if for all
objects $X, Y \in \mathbf{C}$ the morphism $\mathbf{C}(X, Y) \to
\mathbf{D}(FX, FY)$ is a weak equivalence in $\mathbf{V}$. We denote
the class of morphisms in $\Cat(\mathbf{V})$ that are weakly fully
faithful and given by bijections on sets of objects by $\txt{FFB}$.
\end{defn}
The map $\eta_{X} \colon \Cat_{X}(\mathbf{V})[\txt{FF}_{X}^{-1}] \to
\Alg_{\simp^{\op}_{X}}(\mathbf{V}[W^{-1}])$ is natural in $X$, so it
induces a natural transformation of functors $\Set \to
\sSet^{+}$. Applying Corollary~\ref{cor:relcat} we therefore get the
required comparison of ``pre-localized'' homotopy theories:
\begin{thm}\label{thm:AlgcCatCmp}
  The natural transformation $\eta$ induces a functor
  \[\Cat(\mathbf{V})[\txt{FFB}^{-1}] \to
  \AlgCat(\mathbf{V}[W^{-1}])_{\Set}\] and this is an
  equivalence.
\end{thm}

\begin{remark}
  Using Proposition~\ref{propn:totalspmodstr} we can combine the
  (fibrewise) model structures on $\Cat_{X}(\mathbf{V})$ to get a
  model structure on $\Cat(\mathbf{V})$. Explicitly, if $\mathbf{V}$
  is a left proper tractable biclosed monoidal model category
  satisfying the monoid axiom, then there is a model structure on
  $\Cat(\mathbf{V})$ such that a morphism $F \colon \mathbf{C} \to
  \mathbf{D}$ is a weak equivalence \IFF{} $F$ is weakly fully
  faithful and a bijection on objects, and a fibration \IFF{}
  $\mathbf{C}(x,y) \to \mathbf{D}(Fx,Fy)$ is a fibration in
  $\mathbf{V}$ for all $x, y \in \ob \mathbf{C}$. Thus
  $\Cat(\mathbf{V})[\txt{FFB}^{-1}]$ is the \icat{} associated to this
  model category.
\end{remark}

The weakly fully faithful functors that are bijective on objects are
clearly not the right weak equivalences between
$\mathbf{V}$-categories --- just as for ordinary categories the
equivalences are the functors that are fully faithful and essentially
surjective, here they should be the functors that are weakly fully
faithful and essentially surjective up to homotopy, in the following
sense:
\begin{defn}
  Let $\mathbf{V}$ be a monoidal model category. Then the projection
  $\mathbf{V} \to \mathrm{h}\mathbf{V}$ to the homotopy category is a
  monoidal functor; this therefore induces a functor $\Cat(\mathbf{V})
  \to \Cat(\mathrm{h}\mathbf{V})$. A functor of
  $\mathbf{V}$-categories is \emph{homotopically essentially
    surjective} if its image in $\Cat(\mathrm{h}\mathbf{V})$ is
  essentially surjective, and a \emph{DK-equivalence} if it is weakly
  fully faithful and homotopically essentially surjective (or
  equivalently if it induces an equivalence of
  $\mathrm{h}\mathbf{V}$-categories). We write $\txt{DK}$ for the
  class of DK-equivalences in $\Cat(\mathbf{V})$.
\end{defn}

The DK-equivalences in $\Cat(\mathbf{V})$ clearly correspond to the
fully faithful and essentially surjective functors in
$\AlgCat(\mathbf{V}[W^{-1}])_{\Set}$, as defined in \cite[\S
5.2]{enr}. Theorem~\ref{thm:AlgcCatCmp} therefore immediately implies
the following:
\begin{cor}\label{cor:AlgCatDKLoc}
  Suppose $\mathbf{V}$ is a left proper tractable biclosed monoidal
  model category satisfying the monoid axiom. Then
  $\Cat(\mathbf{V})[\txt{DK}^{-1}]$ is equivalent to the
  localization of
  $\AlgCat(\mathbf{V}[W^{-1}])_{\Set}$ with respect to the
  fully faithful and essentially surjective functors.
\end{cor}

Combining this with \cite[Theorem 5.2.17]{enr} we get our main result:
\begin{thm}\label{thm:enrcomp} 
  Suppose $\mathbf{V}$ is a left proper tractable biclosed monoidal
  model category satisfying the monoid axiom. The functor $\eta \colon
  \Cat(\mathbf{V})[\txt{FFB}^{-1}] \to
  \AlgCat(\mathbf{V}[W^{-1}])_{\Set}$ induces an
  equivalence  \[\Cat(\mathbf{V})[\txt{DK}^{-1}] \isoto
  \CatI^{\mathbf{V}[W^{-1}]}.\]
\end{thm}
\begin{proof}
  By \cite[Theorem 5.2.17]{enr}, for any monoidal \icat{}
  $\mathcal{V}$ the localization of $\AlgCat(\mathcal{V})_{\Set}$ at
  the fully faithful and essentially surjective functors is equivalent
  to the corresponding localization of $\AlgCat(\mathcal{V})$, which
  is $\CatIV$ by \cite[Theorem 5.5.6]{enr}. The result follows by
  combining this, in the case where $\mathcal{V}$ is
  $\mathbf{V}[W^{-1}]$, with Corollary~\ref{cor:AlgCatDKLoc}.
\end{proof}

\begin{remark}
  Under the hypotheses of Theorem~\ref{thm:enrcomp} there is a model
  structure on the category $\Cat(\mathbf{V})$ whose weak equivalences
  are the DK-equivalences --- the construction of Muro~\cite{MuroEnr}
  requires slightly weaker hypotheses on $\mathbf{V}$ than our
  theorem. Thus we have shown that $\CatI^{\mathbf{V}[W^{-1}]}$ is the
  \icat{} associated to this model category. Other general
  constructions of model structures on enriched categories are given
  in \cite{HTT,BergerMoerdijkEnr,StanculescuEnr} (see \cite[\S
  1]{BergerMoerdijkEnr} for a historical discussion).
\end{remark}

\begin{ex}
  The \emph{stable model structure} on the category $\Sp^{\Sigma}$ of
  \emph{symmetric spectra}, as described in \cite{HoveyShipleySmith},
  satisfies the hypotheses of Theorem~\ref{thm:enrcomp}. The
  associated monoidal \icat{} is the \icat{} of spectra with the smash
  product monoidal structure. Thus we have an equivalence
  \[ \Cat(\Sp^{\Sigma})[\txt{DK}^{-1}] \isoto \CatI^{\Sp} \]
  between spectral categories and spectral \icats{}.
\end{ex}

\begin{ex}
  The projective model structure on the category $\txt{Ch}^{\geq
    0}(\txt{Mod}_{R})$ of non-negatively graded chain complexes of modules
  over a commutative ring $R$, as described for example in
  \cite{DwyerSpalinski}, satisfies the hypotheses of
  Theorem~\ref{thm:enrcomp}. The same is true of the projective model
  structure on the category $\txt{Ch}(\txt{Mod}_{R})$ of unbounded chain
  complexes of $R$-modules described in \cite[\S
  2.3]{HoveyModCats}. The associated monoidal \icats{} are the bounded
  and unbounded derived \icats{} $\mathcal{D}^{\geq
    0}_{\infty}(\txt{Mod}_{R})$ and $\mathcal{D}_{\infty}(\txt{Mod}_{R})$ of
  $R$-modules, as described in \cite[\S 1.3.2]{HA}. (These are
  equivalent to the \icats{} $\txt{Mod}_{\mathrm{H}R}^{\geq 0}$ and
  $\txt{Mod}_{\mathrm{H}R}$ of connective modules and all modules over the
  Eilenberg-MacLane ring spectrum $\mathrm{H}R$, respectively.) Thus
  we have equivalences
  \[ \Cat(\txt{Ch}^{\geq 0}(\txt{Mod}_{R}))[\txt{DK}^{-1}] \isoto
  \CatI^{\mathcal{D}^{\geq 0}_{\infty}(\txt{Mod}_{R})} \simeq
  \CatI^{\txt{Mod}_{\mathrm{H}R}^{\geq 0}}, \]
  \[ \Cat(\txt{Ch}(\txt{Mod}_{R}))[\txt{DK}^{-1}] \isoto
  \CatI^{\mathcal{D}_{\infty}(\txt{Mod}_{R})} \simeq
  \CatI^{\txt{Mod}_{\mathrm{H}R}}, \] between \icats{} of (two versions
  of) dg-categories and the appropriate corresponding enriched
  \icats{}.
\end{ex}

\section{Comparison with Segal Categories}\label{sec:segcatcomp}
\emph{Segal categories} are a model for the theory of
$(\infty,1)$-categories where composition is only associative up to
coherent homotopy, inspired by Segal's model of
$A_{\infty}$-spaces. They first appeared in papers of Schwänzl and
Vogt~\cite{SchwanzlVogtLoc} and Dwyer, Kan, and
Smith~\cite{DwyerKanSmithDiag}, though not with this name; they were
later rediscovered by Hirschowitz and
Simpson~\cite{HirschowitzSimpson}, who used them as a model for
$(\infty,n)$-categories. A generalization to Segal categories enriched
in a Cartesian model category (i.e. a monoidal model category where
the tensor product is the Cartesian product) was first given by
Pellissier~\cite{Pellissier}, further developed by 
Lurie~\cite{LurieGoodwillie}, and finally extensively studied by
Simpson~\cite{SimpsonSegCats}. In this section we will show that, for
$\mathbf{V}$ a nice Cartesian model category with weak equivalences $W$,
the homotopy theory of Segal categories enriched in $\mathbf{V}$ is
equivalent to that of \icats{} enriched in $\mathbf{V}[W^{-1}]$.
We will first carry out the comparison in the case of a fixed set of
objects, and then apply the results of \S\ref{sec:loc} to prove the
general comparison. 

\begin{defn}
  A model category is \emph{Cartesian} if it is a monoidal model
  category with respect to the Cartesian product. If $\mathbf{V}$ is
  a Cartesian model category, a \emph{$\mathbf{V}$-enriched Segal
    category} (or \emph{Segal $\mathbf{V}$-category}) with set of
  objects $S$ is a functor $\mathbf{C} \colon \simp^{\op}_{S} \to
  \mathbf{V}$ such that for every object $(x_{0}, \ldots, x_{n})$ of
  $\simp^{\op}_{S}$ the \emph{Segal morphism} $\mathbf{C}(x_{0},
  \ldots, x_{n}) \to \mathbf{C}(x_{0},x_{1}) \times \cdots
  \mathbf{C}(x_{n-1},x_{n})$ induced by the projections $(x_{0},
  \ldots, x_{n}) \to (x_{i}, x_{i+1})$ is a weak equivalence. We
  say the Segal category $\mathbf{C}$ is \emph{fibrant} if the objects
  $\mathbf{C}(x_{0},\ldots,x_{n})$ in $\mathbf{V}$ are fibrant for all
  $x_{0},\ldots,x_{n} \in S$, and \emph{strictly unital} if the
  objects $\mathbf{C}(x)$ are final objects in $\mathbf{V}$ for all $x \in S$.
\end{defn}
\begin{remark}
  We can regard $\mathbf{V}$-categories as those Segal categories
  where the Segal morphisms are \emph{isomorphisms}, rather than just
  weak equivalences.
\end{remark}

We can describe fibrant Segal categories with a fixed set $S$ of objects as
the fibrant objects in a Bousfield localization of the projective
model structure on $\Fun(\simp^{\op}_{S}, \mathbf{V})$:
\begin{defn}
  If $X$ is an object of $\simp^{\op}_{S}$, let $i_{X} \colon * \to
  \simp^{\op}_{S}$ denote the functor with image $X$, write $i_{X}^{*}
  \colon \Fun(\simp^{\op}_{S}, \mathbf{V}) \to \mathbf{V}$ for the
  functor given by composition with $i_{X}$, and let $i_{X,!} \colon
  \mathbf{V} \to \Fun(\simp^{\op}_{S}, \mathbf{V})$ be its left
  adjoint, given by left Kan extension along $i_{X}$. Then $i_{X,!}$
  is a left Quillen functor with respect to the projective model
  structure on $\Fun(\simp^{\op}_{S}, \mathbf{V})$. A functor
  $\mathbf{C} \colon \simp^{\op}_{S} \to \mathbf{V}$ is a fibrant
  Segal category \IFF{} it is projectively fibrant and local with
  respect to the morphisms $i_{(x_{0},x_{1}),!}A \coprod \cdots
  \coprod i_{(x_{n-1},x_{n}),!}A \to i_{(x_{0},\ldots,x_{n}),!}A$ for
  all $x_{0},\ldots,x_{n}$ in $S$ and all $A$ in a set of objects that
  generates $\mathbf{V}$ under colimits. If $\mathbf{V}$ is a left
  proper combinatorial Cartesian model category, then we can define a
  model structure whose fibrant objects are fibrant Segal categories
  by taking the left Bousfield localization of the projective model
  structure on $\Fun(\simp^{\op}_{S}, \mathbf{V})$ with respect to
  these morphisms --- this exists under these hypotheses on
  $\mathbf{V}$ by a theorem of Smith (a proof can be found in
  \cite[Theorem 4.7]{BarwickBous}). We refer to this model structure
  as the \defterm{Segal category model structure on functors} and
  write $\Fun(\simp^{\op}_{S}, \mathbf{V})_{\Seg}$ for the category
  $\Fun(\simp^{\op}_{S}, \mathbf{V})$ equipped with this model
  structure.
\end{defn}

To obtain a well-behaved model structure, it turns out to be better to consider
only strictly unital Segal categories. This leads to considering the
category of \emph{$\mathbf{V}$-precategories}:
\begin{defn}
  Let $\mathbf{V}$ be a left proper combinatorial Cartesian model
  category. A \defterm{$\mathbf{V}$-precategory} with set of objects
  $S$ is a functor $\mathbf{C} \colon \simp^{\op}_{S} \to \mathbf{V}$
  such that $\mathbf{C}(x)$ is a final object for all $x \in S$. Write
  $\txt{Precat}_{S}(\mathbf{V})$ for the full subcategory of
  $\Fun(\simp^{\op}_{S}, \mathbf{V})$ spanned by the
  $\mathbf{V}$-precategories  and $u^{*} \colon
  \txt{Precat}_{S}(\mathbf{V}) \to \Fun(\simp^{\op}_{S}, \mathbf{V})$
  for the inclusion. Then $u^{*}$ has a left adjoint, which we denote
  $u_{!}$.
\end{defn}

There is a model structure on $\txt{Precat}_{S}(\mathbf{V})$ analogous
to that for $\Fun(\simp^{\op}_{S}, \mathbf{V})$ we described above:
\begin{propn}[Simpson~{\cite[Propostion 13.4.3]{SimpsonSegCats}}]
  Suppose $\mathbf{V}$ is a left proper combinatorial Cartesian model
  category. There exists a (projective) model structure on
  $\txt{Precat}_{S}(\mathbf{V})$ where a morphism is a weak
  equivalence or fibration if it levelwise is one in $\mathbf{V}$. The
  functor $u^{*} \colon \txt{Precat}_{S}(\mathbf{V}) \to
  \Fun(\simp^{\op}_{S}, \mathbf{V})$ is a right Quillen functor.
\end{propn}

\begin{defn}
  Suppose $\mathbf{V}$ is a left proper combinatorial Cartesian model
  category. The \defterm{(projective) Segal category model structure
    on precategories} is the left Bousfield localization of this
  (projective) model structure on $\txt{Precat}_{S}(\mathbf{V})$
  with respect to the morphisms $u_{!}(i_{(x_{0},x_{1}),!}A \coprod
  \cdots \coprod i_{(x_{n-1},x_{n}),!A}) \to
  u_{!}i_{(x_{0},\ldots,x_{n}),!}A$ for all $(x_{0},\ldots,x_{n})$ in
  $S$ and all $A$ in a set of objects that generates $\mathbf{V}$
  under colimits. We write $\txt{Precat}_{S}(\mathbf{V})_{\Seg}$ for
  the category $\txt{Precat}_{S}(\mathbf{V})$ equipped with this model
  structure.
\end{defn}


Under mild hypotheses these two model categories in the fixed-objects case
are equivalent:
\begin{propn}\label{propn:SegUnitEq}
  Suppose $\mathbf{V}$ is a left proper combinatorial Cartesian model
  category where monomorphisms are cofibrations. Then the adjunction
  $u_{!} \dashv u^{*}$ gives a Quillen equivalence
  \[ \Fun(\simp^{\op}_{S}, \mathbf{V})_{\Seg} \rightleftarrows
  \txt{Precat}_{S}(\mathbf{V})_{\Seg}.\]
\end{propn}
\begin{proof}
  It is obvious that $u^{*}$ is a right Quillen functor, so this is a
  Quillen adjunction. Since $u^{*}$ is fully faithful, the counit
  $u_{!}u^{*}F \to F$ is an isomorphism in
  $\txt{Precat}_{S}(\mathbf{V})$ for all $F$.  By \cite[Lemma
  14.2.1]{SimpsonSegCats} the functor $u_{!}$ only changes the values
  of a functor at the constant sequences $(x, \ldots, x)$ for $x \in
  S$, in which case $u_{!}F$ is given by forming the pushout
  \csquare{F(x)}{*}{F(x,\ldots,x)}{u_{!}F(x,\ldots,x),}{}{F(\sigma)}{}{}
  where $\sigma \colon (x) \to (x,\ldots,x)$ is the map over the
  unique map $s \colon [0] \to [n]$ in $\simp^{\op}$. If $d$ is any
  map $[n] \to [0]$ in $\simp^{\op}$, then $ds = \id$, hence
  $F(\sigma)$ is a monomorphism. By assumption it is therefore a
  cofibration, and so as $\mathbf{V}$ is left proper, the map
  $F(x,\ldots,x) \to u_{!}F(x,\ldots,x)$ is a weak equivalence if
  $F(x) \to *$ is a weak equivalence. Thus $F \to u^{*}u_{!}F$ is a
  levelwise weak equivalence if the map $F(x) \to *$ is a weak equivalence in
  $\mathbf{V}$ for every $x \in S$. Since every object of
  $\Fun(\simp^{\op}_{S}, \mathbf{V})_{\Seg}$ is weakly equivalent to
  one for which this is true, it is clear that the Quillen adjunction $u_{!}
  \dashv u^{*}$ gives an equivalence of homotopy categories, and so is
  a Quillen equivalence.
\end{proof}

Next, we will compare the \icat{} associated to $\Fun(\simp^{\op}_{S},
\mathbf{V})_{\Seg}$ to $\Alg_{\simp^{\op}_{S}}(\mathbf{V}[W^{-1}])$.
We know that the \icat{} associated to the projective model structure
on $\Fun(\simp^{\op}_{S}, \mathbf{V})$ is equivalent to the \icatl{}
functor category $\Fun(\simp^{\op}_{S}, \mathbf{V}[W^{-1}])$. The
Bousfield-localized model category $\Fun(\simp^{\op}_{S},
\mathbf{V})_{\Seg}$ can therefore be identified with the full
subcategory of $\Fun(\simp^{\op}_{S}, \mathbf{V}[W^{-1}])$ spanned by
the objects that are local with respect to certain maps. We can
identify this with the \icat{} of $\simp^{\op}_{S}$-monoids:
\begin{defn}
  Recall that if $\mathcal{V}$ is an \icat{} with finite limits and
  $\mathcal{M}$ is a generalized non-symmetric \iopd{}, an
  \emph{$\mathcal{M}$-monoid} in $\mathcal{V}$ is a functor
  $\mathcal{M} \to \mathcal{V}$ such that for every object $m \in
  \mathcal{M}_{[n]}$, if $m \to m_{i}$ ($i = 1, \ldots, n$) are
  coCartesian morphisms corresponding to the inert maps $\rho_{i}
  \colon [1] \to [n]$ in $\simp$, then the induced morphism $F(m) \to
  F(m_{1}) \times \cdots \times F(m_{n})$ is an equivalence. We write
  $\txt{Mon}_{\mathcal{M}}(\mathcal{V})$ for the full subcategory of
  $\Fun(\mathcal{M}, \mathcal{V})$ spanned by the monoids. There is a
  natural equivalence $\txt{Mon}_{\mathcal{M}}(\mathcal{V}) \simeq
  \Alg_{\mathcal{M}}(\mathcal{V})$ (by \cite[Proposition 3.5.3]{enr}).
\end{defn}

\begin{defn}
  Suppose $\mathcal{V}$ is a presentable \icat{} and $\mathcal{M}$ is
  a generalized non-symmetric \iopd{}. For $m \in \mathcal{M}$, write
  $i_{m} \colon * \to \mathcal{M}$ for the inclusion of this object,
  and let $i_{m,!}$ denote left Kan extension along $i_{m}$. Then for
  any functor $F \colon \mathcal{M} \to \mathcal{V}$ and $X \in
  \mathcal{V}$ we have $\Map(i_{m,!}c_{X}, F) \simeq \Map(c_{X},
  i_{m}^{*}F) \simeq \Map_{\mathcal{V}}(X, F(m))$, where $c_{X}$ is
  the functor $* \to \mathcal{V}$ with image $X$.
\end{defn}
\begin{lemma}
  Suppose $\mathcal{V}$ is a presentable \icat{} such that the
  Cartesian product preserves colimits separately in each variable,
  and $\mathcal{M}$ is a small generalized non-symmetric \iopd{}.
  Then the \icat{} $\txt{Mon}_{\mathcal{M}}(\mathcal{V})$
  is the localization of $\Fun(\mathcal{M}, \mathcal{V})$ with respect
  to the morphisms $i_{m_{1},!}X \amalg \cdots \amalg i_{m_{n},!}X \to
  i_{m,!}X$ for all $m \in \mathcal{M}$ with $X$ ranging over a set of
  objects that generates $\mathcal{V}$ under colimits.
\end{lemma}
\begin{proof}
  A functor $F \colon \mathcal{M} \to \mathcal{V}$ is a monoid \IFF{}
  it is local with respect to these morphisms.
\end{proof}

Since $\txt{Mon}_{\mathcal{M}}(\mathcal{V})$ is
equivalent to $\Alg_{\mathcal{M}}(\mathcal{V})$, we have
proved the following:
\begin{propn}\label{propn:SegComp}
  Suppose $\mathbf{V}$ is a left proper combinatorial Cartesian model
  category, and let $W_{\Seg,S}$ denote the class of weak equivalences
  in $\Fun(\simp^{\op}_{S}, \mathbf{V})_{\Seg}$. Then the natural map
  $\alpha_{S} \colon \Fun(\simp^{\op}_{S},
  \mathbf{V})[W_{\Seg,S}^{-1}] \to
  \Alg_{\simp^{\op}_{S}}(\mathbf{V}[W^{-1}])$ is an equivalence. If
  moreover monomorphisms in $\mathbf{V}$ are cofibrations, then we
  also have a natural equivalence
  $\txt{Precat}_{S}(\mathbf{V})[W_{\txt{Pre},X}^{-1}] \to
  \Alg_{\simp^{\op}_{S}}(\mathbf{V}[W^{-1}])$, where $W_{\txt{Pre},X}$
  denotes the class of weak equivalences in
  $\txt{Precat}_{S}(\mathbf{V})_{\Seg}$.
\end{propn}

Having dealt with the fixed-objects case, we will now allow the set of
objects to vary:
\begin{defn}
  Let $\Seg_{\Fun}(\mathbf{V})$ denote the total space of the right
  Quillen presheaf given by $S \mapsto \Fun(\simp^{\op}_{S},
  \mathbf{V})_{\Seg}$ and let $\txt{Precat}(\mathbf{V})$ denote the
  total space of the right Quillen presheaf given by $S \mapsto
  \txt{Precat}_{S}(\mathbf{V})_{\Seg}$. The adjunction $u_{!} \dashv
  u^{*}$ is natural and so gives a natural transformation
  between these right Quillen presheaves.
\end{defn}

\begin{propn}
  Let $\mathbf{V}$ be a left proper combinatorial Cartesian
  model category. There exist combinatorial model structures on the
  categories $\Seg_{\Fun}(\mathbf{V})$ and $\txt{Precat}(\mathbf{V})$
  where a morphism $F \colon \mathbf{C} \to \mathbf{D}$ is a weak
  equivalence \IFF{} the induced morphism $f$ on objects is a
  bijection and $\mathbf{C} \to f^{*}\mathbf{D}$ is a weak equivalence
  in $\Fun(\simp^{\op}_{\ob \mathbf{C}}, \mathbf{V})_{\Seg}$ or
  $\txt{Precat}_{\ob \mathbf{C}}(\mathbf{V})_{\Seg}$ and a fibration
  \IFF{} $\mathbf{C} \to f^{*}\mathbf{D}$ is a fibration in
  $\Fun(\simp^{\op}_{\ob \mathbf{C}}, \mathbf{V})_{\Seg}$ or
  $\txt{Precat}_{\ob \mathbf{C}}(\mathbf{V})_{\Seg}$. The adjunction
  \[ u_{!} \colon \Seg_{\Fun}(\mathbf{V}) \leftrightarrows
  \txt{Precat}(\mathbf{V}) : u^{*} \] induced by the natural
  transformations $u_{!}$ and $u^{*}$ is a Quillen equivalence.
\end{propn}
\begin{proof}
  This is immediate from Proposition~\ref{propn:totalspmodstr}.
\end{proof}

Now combining Corollary~\ref{cor:relcat} and 
Proposition~\ref{propn:SegComp} we get the following
comparison of ``algebraic'' homotopy theories:
\begin{thm}\label{thm:AlgcSegCmp}
  Suppose $\mathbf{V}$ is a left proper  combinatorial
  Cartesian model category. The natural transformation $\alpha$
  induces a functor $\Seg_{\Fun}(\mathbf{V})[W_{\Fun}^{-1}] \to
  \AlgCat(\mathbf{V}[W^{-1}])_{\Set}$ and this is an
  equivalence, where $W_{\Fun}$ denotes the weak equivalences in the
  model structure on $\Seg_{\Fun}(\mathbf{V})$. If moreover
  monomorphisms in $\mathbf{V}$ are cofibrations, then we also have an
  equivalence $\txt{Precat}(\mathbf{V})[W_{\txt{Precat}}^{-1}] \simeq
  \AlgCat(\mathbf{V}[W^{-1}])_{\Set}$.
\end{thm}

The weak equivalences in $\Seg_{\Fun}(\mathbf{V})$ are difficult to
describe in general; however, a morphism $f \colon \mathbf{C} \to
\mathbf{D}$ between fibrant objects, i.e. Segal categories, is a weak
equivalence \IFF{} it is bijective on objects and a levelwise weak
equivalence --- in fact, given the Segal conditions, it suffices for $f$ to
give a weak equivalence $\mathbf{C}(x,y) \to \mathbf{D}(fx,fy)$ for
all objects $x, y$ in $\mathbf{C}$. To obtain the correct homotopy
theory we clearly also need to invert the morphisms that are fully
faithful and essentially surjective in the appropriate sense:
\begin{defn}
  Composition with the projection $\mathbf{V} \to
  \mathrm{h}\mathbf{V}$ induces a functor \[\Seg_{\Fun}(\mathbf{V}) \to
  \Seg_{\Fun}(\mathrm{h}\mathbf{V}).\] This takes Segal categories to
  categories enriched in $\mathrm{h}\mathbf{V}$. We say a morphism
  between Segal categories in $\Seg_{\Fun}(\mathbf{V})$ is
  \defterm{weakly fully faithful and homotopically essentially
    surjective} if its image in $\Seg_{\Fun}(\mathrm{h}\mathbf{V})$
  corresponds to a fully faithful and essentially surjective functor
  of $\mathrm{h}\mathbf{V}$-categories.
\end{defn}
This definition extends to give a notion of weak equivalence in
$\Seg_{\Fun}(\mathbf{V})$, and similarly in
$\txt{Precat}(\mathbf{V})$; we will refer to these as \emph{Segal
  equivalences}, and denote the class of them as $\txt{SE}$ (in both
$\Seg_{\Fun}(\mathbf{V})$ and $\txt{Precat}(\mathbf{V})$). There are three model structures on
$\txt{Precat}(\mathbf{V})$ with the Segal equivalences as weak
equivalences, namely the \emph{projective}, \emph{injective},
and \emph{Reedy} model structures, constructed in
\cite{SimpsonSegCats}.

The Segal equivalences between Segal categories clearly correspond to the fully
faithful and essentially surjective functors between categorical
algebras, so we get the following:
\begin{propn}
  Suppose $\mathbf{V}$ is a left proper combinatorial
  Cartesian model category. Then there is an equivalence
  \[\Seg_{\Fun}(\mathbf{V})[\txt{SE}^{-1}] \isoto
  \AlgCat(\mathbf{V}[W^{-1}])_{\Set}[\txt{FFES}^{-1}].\]
  If moreover monomorphisms in $\mathbf{V}$ are cofibrations, then
  there is an equivalence
  \[\txt{Precat}(\mathbf{V})[\txt{SE}^{-1}] \isoto
  \AlgCat(\mathbf{V}[W^{-1}])_{\Set}[\txt{FFES}^{-1}].\]
\end{propn}
Combining this with \cite[Theorem 5.2.17]{enr} gives our
comparison result:
\begin{thm}\label{thm:SegCatComp} 
  Suppose $\mathbf{V}$ is a left proper combinatorial Cartesian
  model category. There is an equivalence of \icats{}
  \[ \Seg_{\Fun}(\mathbf{V})[\txt{SE}^{-1}] \isoto
  \CatI^{\mathbf{V}[W^{-1}]}.\]
  If moreover monomorphisms in $\mathbf{V}$ are cofibrations, then
  there is an equivalence
  \[\txt{Precat}(\mathbf{V})[\txt{SE}^{-1}] \isoto
  \CatI^{\mathbf{V}[W^{-1}]}.\]
\end{thm}

\begin{cor}\label{cor:SegnCatComp}
  Let $\mathbf{V}$ be a left proper tractable Cartesian model category
  that is a presheaf category such that the monomorphisms are the
  cofibrations. Then for all $n \geq 0$ there are equivalences of
  \icats{}
  \[ \txt{Precat}^{n}(\mathbf{V})[\txt{SE}^{-1}] \isoto
  \Cat_{(\infty,n)}^{\mathbf{V}[W^{-1}]}.\]
\end{cor}
\begin{proof}
  We wish to apply Theorem~\ref{thm:SegCatComp} inductively. To do
  this we must check that if $\mathbf{V}$ satisfies the given
  hypotheses, then so does a suitable model structure on
  $\txt{Precat}(\mathbf{V})$.
  By \cite[Theorem 21.3.2]{SimpsonSegCats}, if $\mathbf{V}$ is a left
  proper tractable Cartesian  model category then the same is
  true of the \emph{Reedy} model structure on
  $\txt{Precat}(\mathbf{V})$. Moreover, by
  \cite[Proposition 15.7.2]{SimpsonSegCats} if $\mathbf{V}$ is a
  presheaf category such that the monomorphisms are the cofibrations,
  then the injective and Reedy model structures on
  $\txt{Precat}(\mathbf{V})$ coincide, so the Reedy cofibrations are
  the monomorphisms, since these are clearly the injective
  cofibrations. Finally $\txt{Precat}(\mathbf{V})$ is also a presheaf
  category by \cite[Proposition 12.7.6]{SimpsonSegCats}.

  By induction it therefore follows that the Reedy model structure on
  $\txt{Precat}^{n}(\mathbf{V})$ satisfies the hypotheses of
  Theorem~\ref{thm:SegCatComp} for all $n$. Moreover, since the
  monoidal structures on both $\txt{Precat}(\mathbf{V})$ and
  $\CatI^{\mathbf{V}[W^{-1}]}$ are given by the Cartesian product, the
  equivalence between them is automatically an equivalence of
  symmetric monoidal \icats{}, hence induces an equivalence
  $\CatI^{\txt{Precat}(\mathbf{V})} \isoto
  \Cat_{(\infty,2)}^{\mathbf{V}[W^{-1}]}$, etc. By induction we thus
  get a sequence of equivalences
  \[ \txt{Precat}^{n}(\mathbf{V})[\txt{SE}^{-1}] \simeq
  \CatI^{\txt{Precat}^{n-1}(\mathbf{V})[\txt{SE}^{-1}]} \simeq
  \Cat_{(\infty,2)}^{\txt{Precat}^{n-2}(\mathbf{V})[\txt{SE}^{-1}]}
  \simeq \cdots \simeq
  \Cat_{(\infty,n)}^{\mathbf{V}[W^{-1}]}.\qedhere\]
\end{proof}

\begin{ex}
  If we take $\mathbf{V}$ to
  be the category $\sSet$ of simplicial sets, with the usual model
  structure, we get an equivalence
    \[ \txt{Precat}^{n}(\sSet)[\txt{SE}^{-1}] \isoto
  \Cat_{(\infty,n)}, \] where the left-hand side is the \icat{} of the
  $(\infty,n)$-categories of Pellissier-Hirschowitz-Simpson and the
  right-hand side is the \icat{} of $(\infty,n)$-categories defined by
  iterated \icatl{} enrichment.
\end{ex}

\begin{ex}
  We would like to take $\mathbf{V}$ to be the category $\Set$ of
  sets, equipped with the trivial model structure, but of course this
  does not satisfy the hypothesis that cofibrations are
  monomorphisms. We therefore need to consider instead a model
  category $\mathbf{M}$, Quillen equivalent to $\Set$, that does
  satisfy the hypotheses of the theorem. For example, following
  \cite[\S 22.1]{SimpsonSegCats} we can let $\mathbf{M}$ be an
  appropriate localization of the Reedy model structure on
  $\txt{Precat}^{2}(*)$, or we can take $\mathbf{M}$ to be the
  Bousfield localization of the usual model structure on $\sSet$ with
  respect to the morphisms $\partial \Delta^{n} \to \Delta^{0}$ for
  all $n \geq 2$. We then get an equivalence
  \[ \txt{Precat}^{n}(\Set)[\txt{SE}^{-1}] \isoto
  \txt{Precat}^{n}(\mathbf{M})[\txt{SE}^{-1}] \isoto
  \Cat_{(\infty,n)}^{\mathbf{M}[W^{-1}]} \isoto \Cat_{n}, \] where the
  left-hand side is the \icat{} of Tamsamani's
  $n$-categories~\cite{Tamsamani} and the right-hand side is the
  \icat{} of $n$-categories defined by iterated \icatl{} enrichment.
\end{ex}

\section{Comparison with Iterated Segal Spaces}\label{sec:compitsegsp}
We saw in the previous section that the \icat{} $\Cat_{(\infty,n)}$ of
$(\infty,n)$-categories, obtained by iterated enrichment, is
equivalent to that associated to the model category of $n$-fold Segal
categories, which is another model for the homotopy theory of
\incats{}. Since this model is known to satisfy the axioms of Barwick
and Schommer-Pries~\cite{BarwickSchommerPriesUnicity}, it follows that
$\Cat_{(\infty,n)}$ is equivalent to all the usual models for
\incats{}. However, this comparison was somewhat indirect. Our goal in
this section is to give a more direct comparison between
$\Cat_{(\infty,n)}$ and another established model of
$(\infty,n)$-categories, namely the iterated Segal spaces of Barwick~\cite{BarwickThesis}.

We will deduce this comparison from a slightly more general result: we
will prove that if $\mathcal{X}$ is an \emph{absolute distributor}, in
the sense of \cite{LurieGoodwillie}, then categorical algebras in
$\mathcal{X}$ are equivalent to Segal spaces in $\mathcal{X}$, and
complete categorical algebras are equivalent to complete Segal
spaces. We begin with a brief review of the notion of distributor:
\begin{defn}
  A \emph{distributor} consists of an \icat{} $\mathcal{X}$ together
  with a full subcategory $\mathcal{Y}$ such that:
  \begin{enumerate}[(1)]
  \item The \icats{} $\mathcal{X}$ and $\mathcal{Y}$ are presentable.
  \item The full subcategory $\mathcal{Y}$ is closed under small
    limits and colimits in $\mathcal{X}$.
  \item If $X \to Y$ is a morphism in $\mathcal{X}$ such that $Y \in
    \mathcal{Y}$, then the pullback functor $\mathcal{Y}_{/Y} \to
    \mathcal{X}_{/X}$ preserves colimits.
  \item Let $\mathcal{O}$ denote the full subcategory of
    $\Fun(\Delta^{1}, \mathcal{X})$ spanned by those morphisms $f
    \colon X \to Y$ such that $Y \in  \mathcal{Y}$, and let $\pi
    \colon \mathcal{O} \to \mathcal{Y}$ be the functor given by
    evaluation at $1 \in \Delta^{1}$. Since $\mathcal{X}$ admits
    pullbacks, the evaluation-at-1 functor $\Fun(\Delta^{1}, \mathcal{X})
    \to \mathcal{X}$ is a Cartesian fibration, hence so is $\pi$. Let
    $\chi \colon \mathcal{Y} \to \LCatI^{\op}$ be a functor that
    classifies $\pi$. Then $\chi$ preserves small limits.
  \end{enumerate}
\end{defn}

\begin{defn}
  An \emph{absolute distributor} is a presentable \icat{}
  $\mathcal{X}$ such that the unique colimit-preserving functor $\mathcal{S}
  \to \mathcal{X}$ that sends $*$ to the final object is fully
  faithful, and $\mathcal{S} \subseteq \mathcal{X}$ is a distributor.
\end{defn}

Now we can recall the definition of a Segal space in an absolute distributor:
\begin{defn}
  Suppose $\mathcal{C}$ is an \icat{} with finite limits. A
  \defterm{category object} in $\mathcal{C}$ is a simplicial object $F
  \colon \simp^{\op}\to \mathcal{C}$ such that for each $n$ the map 
  \[ F_{n} \to F_{1} \times_{F_{0}} \cdots \times_{F_{0}} F_{1} \]
  induced by the inclusions $\{i,i+1\} \hookrightarrow [n]$ and $\{i\}
  \hookrightarrow [n]$ is an equivalence.
\end{defn}

\begin{defn}
  Let $\mathcal{X}$ be an absolute distributor. A \defterm{Segal
    space} in $\mathcal{X}$ is a category object $F \colon \simp^{\op}
  \to \mathcal{X}$ such that $F([0])$ is in $\mathcal{S} \subseteq
  \mathcal{X}$.
\end{defn}

Our goal is now to prove the following:
\begin{thm}\label{thm:DistCatAlg}
  Suppose $\mathcal{X}$ is an absolute distributor. There is an
  equivalence \[\AlgCat(\mathcal{X}) \isoto \Seg(\mathcal{X}),\] given
  by sending a $\simp^{\op}_{S}$-algebra $\mathcal{C}$ to the left Kan
  extension $\pi_{!}\mathcal{C}'$ of the composite
  \[ \mathcal{C}' \colon \simp^{\op}_{S} \xto{\mathcal{C}}
  \mathcal{X}^{\times} \to \mathcal{X} \] along $\pi \colon
  \simp^{\op}_{S} \to \simp^{\op}$, where the second map (which sends
  $(S_{1},\ldots, S_{n}) \in \mathcal{X}^{\times}_{[n]}$ to $S_{1}
  \times \cdots \times S_{n}$) comes from a Cartesian structure in the
  sense of \cite[Definition 2.4.1.1]{HA}.
\end{thm}

For the proof we need some more technical results:
\begin{propn}[{\cite[Corollary 1.2.5]{LurieGoodwillie}}]\label{propn:DistrProp}
  Suppose $\mathcal{Y} \subseteq \mathcal{X}$ is a distributor. Let
  $K$ be a small simplicial set, and let $\overline{\alpha} \colon
  \overline{p} \to \overline{q}$ be a natural transformation between
  functors $\overline{p}, \overline{q} \colon K^{\triangleright} \to
  \mathcal{X}$. If $\overline{q}$ is a colimit diagram in
  $\mathcal{Y}$ and $\alpha = \overline{\alpha}|_{K}$ is Cartesian,
  then $\overline{\alpha}$ is Cartesian \IFF{} $\overline{p}$ is a
  colimit diagram.
\end{propn}

\begin{lemma}\label{lem:DistrOverSpace}
  Suppose $\mathcal{X}$ is an absolute distributor. Then for
  every space $X \in \mathcal{S}$, the map \[\gamma_{X} \colon \Fun(X,
  \mathcal{X}) \to \mathcal{X}_{/X}\] that sends a functor $F \colon X
  \to \mathcal{X}$ to its colimit is an equivalence of \icats{}.
\end{lemma}
\begin{proof}
  Let $\xi \colon X \to \mathcal{X}$ be the constant functor at the
  final object $* \in \mathcal{S} \subseteq \mathcal{X}$. Since $X$ is
  a space, a functor $F \colon X \to \mathcal{X}$ sends every morphism
  in $X$ to an equivalence in $\mathcal{X}$, and so the unique natural
  transformation $F \to \xi$ is Cartesian.

  Write $\overline{\xi} \colon X^{\triangleright} \to \mathcal{X}$ for a
  colimit diagram extending $\xi$. Then $\gamma_{X}$ factors as
  \[ \Fun(X, \mathcal{X}) \simeq \Fun(X, \mathcal{X})_{/\xi}
  \xto{\phi_{1}} \Fun(X^{\triangleright},
  \mathcal{X})_{/\overline{\xi}} \xto{\phi_{2}} \mathcal{X}_{/X},\]
  where $\phi_{2}$ is given by evaluation at the cone point. The
  functor $\phi_{1}$ gives an equivalence between $\Fun(X,
  \mathcal{X})_{/\xi}$ and the full subcategory $\mathcal{E}_{1}$ of
  $\Fun(X^{\triangleright}, \mathcal{X})_{/\overline{\xi}}$ spanned by
  the colimit diagrams. On the other hand, the restriction of
  $\phi_{2}$ to the full subcategory $\mathcal{E}_{2}$ spanned by the
  Cartesian natural transformations to $\overline{\xi}$ is also
  clearly an equivalence. By Proposition~\ref{propn:DistrProp} the
  subcategories $\mathcal{E}_{1}$ and $\mathcal{E}_{2}$ coincide, and
  so the composite $\gamma_{X}$ is indeed an equivalence.
\end{proof}

\begin{propn}\label{propn:FuntoXEq}
  Let $\mathcal{O}$ be an \icat{}, and let $F \colon \mathcal{O} \to
  \mathcal{S}$ be a functor; write $\pi \colon \mathcal{O}_{F} \to
  \mathcal{O}$ for the left fibration associated to $F$. Suppose
  $\mathcal{X}$ is an absolute distributor. Then left Kan
  extension along $\pi$ gives an equivalence \[\Fun(\mathcal{O}_{F},
  \mathcal{X}) \isoto \Fun(\mathcal{O}, \mathcal{X})_{/F}.\]
\end{propn}
\begin{proof}
  By \cite[Proposition A.1.5]{enr} the \icat{}
  $\Fun(\mathcal{O}_{F}, \mathcal{X})$ is equivalent to the \icat{} of
  sections of the Cartesian fibration $\mathcal{E} \to \mathcal{O}$
  whose fibre at $X \in \mathcal{O}$ is $\Fun(F(X),
  \mathcal{X})$. Since $\mathcal{X}$ is an absolute distributor,
  by Lemma~\ref{lem:DistrOverSpace} the \icat{} $\mathcal{E}$ is
  equivalent over $\mathcal{O}$ to the total space $\mathcal{E}'$ of
  the Cartesian fibration associated to the functor sending $X$ to
  $\mathcal{X}_{/F(X)}$. Then $\mathcal{E}'$ is the pullback along $F$
  of the
  Cartesian fibration $\Fun(\Delta^{1}, \mathcal{X}) \to \mathcal{X}$
  given by evaluation at $1$, so we have an equivalence between the
  \icat{} $\Fun_{\mathcal{O}}(\mathcal{O}, \mathcal{E}')$ of sections
  and the fibre of $\Fun(\mathcal{O} \times \Delta^{1}, \mathcal{X})
  \simeq \Fun(\Delta^{1}, \Fun(\mathcal{O}, \mathcal{X})) \to
  \Fun(\mathcal{O}, \mathcal{X})$ at $F$. This is clearly equivalent
  to $\Fun(\mathcal{O}, \mathcal{X})_{/F}$, which completes the proof.
\end{proof}


\begin{propn}\label{propn:DistrSegXEq}
  Let $S$ be a space, and let $\pi \colon \simp^{\op}_{S} \to
  \simp^{\op}$ be the usual projection. Let $\pi_{!} \colon
  \Fun(\simp^{\op}_{S}, \mathcal{X}) \to \Fun(\simp^{\op},
  \mathcal{X})$ be the functor given by left Kan extension along
  $\pi$. Then a functor $F \colon \simp^{\op}_{S}\to \mathcal{X}$ is a
  $\simp^{\op}_{S}$-monoid \IFF{} $\pi_{!}F$ is a Segal space.
\end{propn}
\begin{proof}
  It is clear that $\pi_{!}F([0])$ is equivalent to $S$. We must thus
  show that the Segal morphism
  \[ \pi_{!}F([n]) \to \pi_{!}F([1]) \times_{S} \cdots \times_{S}
  \pi_{!}F([1]) =: (\pi_{!}F)^{\Seg}_{[n]} \] is an equivalence \IFF{}
  $F$ is a $\simp^{\op}_{S}$-monoid. Since $\pi$ is a coCartesian
  fibration, we have an equivalence $\pi_{!}F([n]) \simeq \colim_{\xi
    \in S^{\times (n+1)}} F(\xi)$. It thus suffices to show that
  $(\pi_{!}F)^{\Seg}_{[n]}$ is also a colimit of this diagram \IFF{}
  $F$ is a $\simp^{\op}_{S}$-monoid. There is a natural transformation
  $(S^{\times (n+1)})^{\triangleright} \to \Fun(\Delta^{1},
  \mathcal{X})$ that sends $\xi \in S^{\times (n+1)}$ to $F(\xi) \to
  \xi$ and $\infty$ to $(\pi_{!}F)^{\Seg}_{[n]} \to S^{\times (n+1)}$.
  Since $\mathcal{X}$ is an absolute distributor, by
  Proposition~\ref{propn:DistrProp} the colimit is $(\pi_{!}F)^{\Seg}_{[n]}$ \IFF{} this
  natural transformation is Cartesian. Since $S^{\times (n+1)}$ is a
  space, this is equivalent to the square
  \nolabelsmallcsquare{F(\xi)}{(\pi_{!}F)^{\Seg}_{[n]}}{\xi}{S^{\times
      (n+1)}} being a pullback square for all $\xi$, so we are reduced to
  showing that the fibre of $(\pi_{!}F)^{\Seg}_{[n]} \to S^{\times
    (n+1)}$ at $\xi$ is $F(\xi)$ \IFF{} $F$ is a
  $\simp^{\op}_{S}$-monoid. Since limits commute, if $\xi$ is $(s_{0},
  \ldots, s_{n})$ this fibre is the iterated fibre product
  \[ (\pi_{!}F[1])_{(s_{0},s_{1})} \times_{(\pi_{!}F[0])_{(s_{1})}}
  \cdots \times_{(\pi_{!}F[0])_{(s_{n-1})}}
  (\pi_{!}F[1])_{(s_{n-1},s_{n})}.\]
  But using Proposition~\ref{propn:DistrProp} again it is clear that the
  natural maps $F(x,y) \to (\pi_{!}F[1])_{(x,y)}$ and $* \simeq F(x)
  \to (\pi_{!}F)_{(x)}$ are equivalences for all $x,y \in S$. Thus the
  map $F(\xi) \to (\pi_{!}F)^{\Seg}_{[n],\xi}$ is equivalent to the
  natural map
  \[ F(\xi) \to F(s_{0},s_{1}) \times \cdots \times F(s_{n-1},
  s_{n}).\] By definition this is an equivalence for all $\xi \in
  \simp^{\op}_{S}$ \IFF{} $F$ is a $\simp^{\op}_{S}$-monoid, which
  completes the proof.
\end{proof}

\begin{defn}
  Let $i \colon * \to \simp^{\op}$ denote the inclusion of the object
  $[0]$. Then composition with $i$ gives a functor $i^{*} \colon
  \Seg(\mathcal{X}) \to \mathcal{S}$ with left and right adjoints
  $i_{!}$ and $i_{*}$, given respectively by left and right Kan
  extension. Observe that by definition $\simp^{\op}_{X} \to
  \simp^{\op}$ is the left fibration associated to $i_{*}X \in
  \Seg(\mathcal{S})$.
\end{defn}

\begin{cor}
  Let $S$ be a space, and let $\pi \colon \simp^{\op}_{S} \to
  \simp^{\op}$ denote the canonical projection. By
  Proposition~\ref{propn:FuntoXEq} the functor
  \[ \pi_{!} \colon \Fun(\simp^{\op}_{S}, \mathcal{X}) \to
  \Fun(\simp^{\op}, \mathcal{X})_{/i_{*}S} \] given by left Kan
  extension is an equivalence. 

  Under this equivalence, the full subcategory
  $\Mon_{\simp^{\op}_{S}}(\mathcal{X})$ of $\simp^{\op}_{S}$-monoids
  corresponds to the full subcategory of $\Fun(\simp^{\op},
  \mathcal{X})_{/i_{*}S}$ spanned by the Segal spaces $Y_{\bullet}$
  such that $Y_{0}\simeq S$ and the map $Y_{\bullet} \to i_{*}S$ is
  given by the adjunction unit $Y_{\bullet} \to i_{*}i^{*}Y_{\bullet}
  \simeq i_{*}S$.
\end{cor}
\begin{proof}
  It is clear that $\pi_{!}$ takes
  $\Mon_{\simp^{\op}_{S}}(\mathcal{X})$ into the full
  subcategory of $\Fun(\simp^{\op}, \mathcal{X})_{/i_{*}S}$ spanned by
  simplicial spaces $Y_{\bullet}$ with $Y_{0} \simeq S$ and the map
  $Y_{\bullet} \to i_{*}S$ given by the adjunction unit $Y_{\bullet}
  \to i_{*}i^{*}Y \simeq i_{*}S$. The result therefore follows by
  Proposition~\ref{propn:DistrSegXEq}.
\end{proof}

\begin{cor}\label{cor:DopSMndEq2} 
  Let $S$ be a space, and let $\pi \colon \simp^{\op}_{S} \to
  \simp^{\op}$ denote the canonical projection. The functor $\pi_{!}
  \colon \Fun(\simp^{\op}_{S}, \mathcal{X}) \to \Fun(\simp^{\op},
  \mathcal{X})$ given by left Kan extension along $\pi$ gives an
  equivalence of the full subcategory
  $\Mon_{\simp^{\op}_{S}}(\mathcal{X})$ of $\simp^{\op}_{S}$-monoids
  with the subcategory $\Seg(\mathcal{X})_{S}$ of Segal spaces with
  $0$th space $S$ and morphisms that are the identity on the $0$th
  space.
\end{cor}

\begin{proof}[Proof of Theorem~\ref{thm:DistCatAlg}]
  If $\mathcal{V}$ is an \icat{} with finite products, pulling
  back the monoid fibration $\Mon(\mathcal{V}) \to \OpdIns$ of
  \cite[Remark 3.6.3]{enr} along $\simp^{\op}_{(\blank)}$ gives a
  Cartesian fibration $\Mon_{\txt{cat}}(\mathcal{V})$ with an
  equivalence \[\AlgCat(\mathcal{V}) \isoto
  \Mon_{\txt{cat}}(\mathcal{V})\] over $\mathcal{S}$. Taking left Kan
  extensions along the projections $\simp^{\op}_{S}\to \simp^{\op}$
  for all $S \in \mathcal{S}$ we get (using
  Proposition~\ref{propn:DistrSegXEq}) a commutative square
  \opctriangle{\Mon_{\txt{cat}}(\mathcal{X})}{\Seg(\mathcal{X})}{\mathcal{S}.}{\Phi}{}{\txt{ev}_{[0]}}
  By \cite[Lemma A.1.6]{enr} it is clear that $\txt{ev}_{[0]} \colon
  \Seg(\mathcal{X}) \to \mathcal{S}$ is a Cartesian fibration, and the functor
  $\Phi$ preserves Cartesian morphisms by Proposition~\ref{propn:DistrProp}.
  It thus suffices to prove that for each $S \in
  \mathcal{S}$ the functor on fibres
  $\Mon_{\simp^{\op}_{S}}(\mathcal{X}) \to \Seg(\mathcal{X})_{S}$ is an
  equivalence, which is the content of Corollary~\ref{cor:DopSMndEq2}.
\end{proof}

Our goal is now to deduce that the equivalence of
Theorem~\ref{thm:DistCatAlg} induces an equivalence between complete
categorical algebras and complete Segal spaces. We will first review
the definition of the latter:
\begin{defn}
  Write $\txt{Gpd}(\mathcal{S})$ for the full subcategory of
  $\Seg(\mathcal{S})$ spanned by the \emph{groupoid objects}, i.e. the
  simplicial objects $X$ such that for every partition $[n] = S \cup
  S'$ where $S \cap S'$ consists of a single element, the diagram
  \nolabelsmallcsquare{X([n])}{X(S)}{X(S')}{X(S \cap S')} is a
  pullback square. Let $\mathcal{X}$ be an absolute distributor, and
  let $\Lambda \colon \mathcal{X} \to \mathcal{S}$ denote the right
  adjoint to the inclusion $\mathcal{S} \hookrightarrow
  \mathcal{X}$. The inclusion $\txt{Gpd}(\mathcal{S}) \hookrightarrow
  \Seg(\mathcal{S}) \hookrightarrow \Seg(\mathcal{X})$ admits a right
  adjoint $\iota \colon \Seg(\mathcal{X}) \to \txt{Gpd}(\mathcal{S})$,
  which is the composite of the functor $\Lambda \colon
  \Seg(\mathcal{X}) \to \Seg(\mathcal{S})$ induced by $\Lambda$, and
  $\iota \colon \Seg(\mathcal{S}) \to \txt{Gpd}(\mathcal{S})$. We say
  a Segal space $F \colon \simp^{\op} \to \mathcal{X}$ is
  \defterm{complete} if the groupoid object $\iota F$ is constant.
\end{defn}

\begin{remark}
  By \cite[Lemma 5.1.14]{enr}, a Segal space $F$ is complete \IFF{} the
  map \[\iota F(s^{0}) \colon \iota F[0] \to \iota F[1]\] is an
  equivalence.
\end{remark}

\begin{defn}
  Let $E^{n}$ denote the Segal space $i_{*}\{0, \ldots, n\}$. If
  $\mathcal{X}$ is an absolute distributor we also write $E^{n}$
  for $E^{n}$ regarded as a Segal space in $\mathcal{X}$ via the
  inclusion $\mathcal{S} \hookrightarrow \mathcal{X}$.
\end{defn}

\begin{propn}\label{propn:SegXcomplloc}
  Suppose $\mathcal{X}$ is an absolute distributor. Then a Segal
  space $F$ in $\mathcal{X}$ is complete \IFF{} it is local with
  respect to the morphism $E^{1}\to E^{0}$.
\end{propn}
\begin{proof}
  It is clear that $F$ is local with respect to $E^{1} \to E^{0}$,
  considered as a morphism in $\Seg(\mathcal{X})$, \IFF{} the Segal
  space $\Lambda F$ in $\mathcal{S}$ is local with respect to $E^{1}
  \to E^{0}$, considered as a morphism in $\Seg(\mathcal{S})$. On the
  other hand, $F$ is complete \IFF{} $\Lambda F$ is complete, so it
  suffices to prove this for Segal spaces in $\mathcal{S}$. This case
  is part of \cite[Proposition 6.4]{RezkCSS}.
\end{proof}

\begin{defn}
  Let $\txt{CSS}(\mathcal{X})$ denote the full subcategory of
  $\Seg(\mathcal{X})$ spanned by the complete Segal spaces; by
  Proposition~\ref{propn:SegXcomplloc} this is the localization of
  $\Seg(\mathcal{X})$ with respect to the morphism $E^{1} \to E^{0}$.
\end{defn}

\begin{thm}\label{thm:CatIeqCSS}
  Let $\mathcal{X}$ be an absolute distributor. The equivalence
  $\AlgCat(\mathcal{X}) \isoto \Seg(\mathcal{X})$ induces an
  equivalence $\CatI^{\mathcal{X}} \isoto \txt{CSS}(\mathcal{X})$.
\end{thm}
\begin{proof}
  It is clear that $E^{n}_{\mathcal{X}} \in
  \AlgCat(\mathcal{X})$ corresponds to $E^{n} \in
  \Seg(\mathcal{X})$ under this equivalence. Both sides are
  therefore the localization with respect to $E^{1} \to E^{0}$.
\end{proof}

\begin{defn}
  By \cite[Corollary 1.3.4]{LurieGoodwillie}, if $\mathcal{X}$ is an
  absolute distributor, then $\txt{CSS}(\mathcal{X})$ is also an
  absolute distributor. We therefore have absolute distributors
  $\txt{CSS}^{n}(\mathcal{X})$ of \emph{$n$-fold complete
    Segal spaces} in $\mathcal{X}$.
\end{defn}

Applying Theorem~\ref{thm:CatIeqCSS} inductively, we get:
\begin{cor}
  Let $\mathcal{X}$ be an absolute distributor. Then
  $\Cat_{(\infty,n)}^{\mathcal{X}} \simeq \txt{CSS}^{n}(\mathcal{X})$.
\end{cor}

In particular, taking $\mathcal{X}$ to be the \icat{} $\mathcal{S}$ of
spaces, we obtain the desired comparison with iterated Segal spaces:
\begin{cor}
  There is an equivalence $\Cat_{(\infty,n)} \simeq
  \txt{CSS}^{n}(\mathcal{S})$.
\end{cor}

\bibliographystyle{gtart}


\end{document}